\documentclass[11pt]{amsart}
\usepackage{amsmath}
\usepackage{amsxtra}
\usepackage{amscd}
\usepackage{amsthm}
\usepackage{amsfonts}
\usepackage{amssymb}
\usepackage{eucal}

\newtheorem{cor}[subsection]{Corollary}
\newtheorem{lem}[subsection]{Lemma}
\newtheorem{prop}[subsection]{Proposition}

\newtheorem{thm}[subsection]{Theorem}
\newtheorem{defn}[subsection]{Definition}
\theoremstyle{remark}

\newcommand{\thmref}[1]{Theorem~\ref{#1}}
\newcommand{\secref}[1]{Section~\ref{#1}}
\newcommand{\lemref}[1]{Lemma~\ref{#1}}
\newcommand{\propref}[1]{Proposition~\ref{#1}}
\newcommand{\corref}[1]{Corollary~\ref{#1}}

\newcommand{\nc}{\newcommand}

\nc{\ssec}{\subsection}
\nc{\sssec}{\subsubsection}

\nc{\on}{\operatorname}

\nc{\ZZ}{{\mathbb Z}}
\nc{\Z}{{\mathcal C}}
\nc{\z}{{\mathbf c}}
\nc{\TT}{{\mathcal T}}
\nc{\TTb}{{\overline{\mathbb T}}}
\renewcommand{\tt}{{\mathbf t}}
\nc{\T}{{\mathcal T}}
\nc{\Tf}{{\mathfrak T}}
\nc{\tof}{{\mathfrak t}}
\nc{\tf}{{\mathfrak t}}

\nc{\NN}{{\mathbb N}}
\nc{\CC}{{\mathbb C}}
\nc{\GG}{{\mathbb G}}
\nc{\DD}{{\mathbb D}}
\nc{\PP}{{\mathbb P}}

\nc{\A}{{\mathcal A}}
\nc{\B}{{\mathcal Fl}}
\nc{\Bu}{{\mathcal B}}
\nc{\R}{{\mathcal R}}
\nc{\N}{{\mathcal N}}
\nc{\C}{{\mathcal Q}}
\renewcommand{\O}{{\mathcal O}}
\nc{\K}{{\mathcal K}}
\nc{\M}{{\mathcal M}}
\nc{\V}{{\mathcal V}}

\renewcommand{\P}{{\mathcal P}}
\renewcommand{\H}{{\mathcal H}}
\nc{\U}{{\mathcal U}}
\renewcommand{\L}{{\mathcal L}}
\nc{\G}{{\mathcal G}}
\nc{\Y}{{\mathcal Y}}

\nc{\Hf}{{\mathfrak H}}
\nc{\Mf}{{\mathfrak M}}
\nc{\Xf}{{\mathfrak X}}
\nc{\gf}{{\mathfrak g}}
\nc{\af}{{\mathfrak a}}
\nc{\qf}{{\mathfrak q}}
\nc{\pf}{{\mathfrak p}}
\nc{\mf}{{\mathfrak m}}
\nc{\bof}{{\mathfrak b}}
\renewcommand\sl{{\mathfrak s}{\mathfrak l}}
\nc\gl{{\mathfrak g}{\mathfrak l}}
\nc\ev{\mathbf{ev}}

\nc{\uf}{{\mathfrak u}}
\nc{\zf}{{\mathfrak z}}

\nc\tilgf{{\widetilde\gf}}
\nc\tili{{\widetilde i}}
\nc\tilL{{\widetilde L}}
\nc\tilE{{\widetilde E}}
\nc\tilU{{\widetilde U}}
\nc\tilS{{\widetilde S}}
\nc\tils{{\widetilde s}}
\nc\tilX{{\widetilde X}}
\nc\tily{{\widetilde y}}
\nc\tilz{{\widetilde z}}
\nc\tilrho{{\widetilde\rho}}
\nc\tilphi{{\widetilde\phi}}
\nc\tilsigma{{\widetilde\sigma}}
\nc\tilpi{{\widetilde\pi}}
\nc\tilZ{{\widetilde \Z}}
\nc\tilB{{\widetilde B}}
\nc\tilY{{\widetilde Y}}
\nc\bfHiggs{{\mathbf H}{\mathbf i}{\mathbf g}{\mathbf g}{\mathbf s}}
\nc\bfHitch{{\mathbf H}{\mathbf i}{\mathbf t}{\mathbf c}{\mathbf h}}
\nc\bfCam{{\mathbf C}{\mathbf a}{\mathbf m}}
\nc\bfCov{{\mathbf C}{\mathbf o}{\mathbf v}}
\nc\bfBun{{\mathbf B}{\mathbf u}{\mathbf n}}
\nc\bfReg{{\mathbf R}{\mathbf e}{\mathbf g}}
\nc\GNb{\overline{G/N}}
\nc\GTb{\overline{G/T}}
\nc\BTb{\overline{B/T}}
\nc\Tb{\overline{T}}
\nc\Xb{\overline{X}}
\newcommand\Sbar{\overline{S}}
\nc\chib{\overline{\chi}}

\nc\Bun{\on{Bun}}

\nc\Spec{\on{Spec}}
\nc\Ind{\on{Ind}}
\nc\Hom{\on{Hom}}
\nc\Mor{\on{Mor}}
\nc\Aut{\on{Aut}}
\nc\et{\on{Sch}_{et}}
\nc\Tors{\on{Tors}}
\nc\one{{\mathbf 1}}
\nc\Pic{\on{Pic}}
\nc\bi{\bibitem}
\nc\Higgs{\on{Higgs}}
\nc\Higgsb{\overline{\Higgs}}
\nc\Cam{\on{Cam}}
\nc\toeq{\stackrel{\sim}{\longrightarrow}}
\title{The gerbe of Higgs bundles}

\author{R.Y.~Donagi and D.~Gaitsgory}

\begin{document}

\maketitle

\tableofcontents
\setcounter{section}{-1}
\section{Introduction}

The purpose of this work is to describe the (category of) Higgs bundles on
a scheme $X / \CC$ having a given cameral cover $\tilX$.  We show that this
category is a $T_{\tilX}$-gerbe, where $T_{\tilX}$ is a certain sheaf of
abelian groups on $X$, and we describe the class of this gerbe precisely.
In particular, it follows that the set of isomorphism classes of Higgs
bundles with a fixed cameral cover $\tilX$ is a torsor over the group
$H^1(X,T_{\tilX})$, which itself parametrizes $T_{\tilX}$-torsors on $X$.
This underlying group $H^1(X,T_{\tilX})$ can be described as a generalized
Prym variety, whose connected component is either an abelian variety or a
degenerate abelian variety.

The hardest part of our work, though, goes into identifying precisely the
$H^1(X,T_{\tilX})$-torsor we get, or in other words, identifying the class
of the gerbe.  This class is surprisingly complicated. One piece of it can
be identified as a twist along the ramification divisors of ${\tilX}$ over
$X$, and is present for all groups $G$. A second piece is a shift which can
be present even for unramified covers. While the twist along the
ramification expresses properties of the cameral cover, this shift
expresses the non-vanishing of a certain group cohomology element -
specifically, the extension class $[N]$ of the normalizer $N=N_G(T)$,
which is an element in the cohomology group $H^2(W,T)$ of the Weyl
group acting on the maximal torus.  It vanishes for some groups, such as
$GL(n), PGL(n), SL(2n+1), SO(n)$, but not for others such as $SL(2n)$.  Yet
a third piece is present only for the groups $SO(2n+1)$ (or groups
containing them as direct factors); this piece expresses the existence of
non-primitive coroots, which amounts to the non vanishing of an element in
another cohomology group.  We give several examples to illustrate these
individual ingredients as well as their combined effect.

Throughout this work, we let $G$ be a connected reductive group, and let
$X$ be a scheme over the complex numbers. A Higgs bundle over $X$ is a
principal $G$-bundle plus some additional data. We describe this
additional data next: first for $G=GL(n)$, and then
for all $G$, in subsection \ref{0.1}.
In the remainder of this introduction we will outline
our results \ref{0.2}, discuss some examples and applications \ref{0.3},
and review some related results in the literature \ref{0.4}. The
notation we employ is summarized in \ref{0.5}.

\ssec{Abelianization: Higgs bundles and cameral covers} \label{0.1}
It is especially easy to spell out the definition when $G=GL(n)$.
In this case a $G$-bundle is the same as a vector bundle $E$ over $X$, and
a Higgs structure on it is a subbundle of commutative associative algebras
$\z_X\subset\on{End}_{\O_X}(E)$, which has rank $n$ over $X$ and such that
$\z_X$ is locally generated by one section. In this case the relative
spectrum
of $\z_X$ over $X$ is a flat $n$-sheeted cover of $X$, called the
{\it spectral cover} corresponding to our Higgs bundle. We will
denote it by $\Xb$.

How can we classify
Higgs bundles with a given spectral cover $\Xb$?
The answer is simple: these are in
bijection with line bundles on $\Xb$. Thus,
by asking not just for principal $G$-bundles,
but rather for $G$-bundles endowed with a Higgs structure with a
fixed spectral cover, we go from
a {\it non-abelian} problem to an {\it abelian} one.

The natural
question  now is how to extend the above discussion to other reductive
groups. It turns out that the notion of an abstract Higgs bundle is
quite easy to generalize. Namely, a Higgs bundle is a pair
$(E_G,\z_X)$, where $E_G$ is a principal $G$-bundle over $X$ and $\z_X$
is a subbundle of the associated bundle of Lie algebras
${\gf}_{E_{G}}$, whose fibers are {\it regular centralizers}. The
precise definition is given in \secref{intrHiggs}. Here we only recall
that a regular centralizer in the Lie algebra ${\gf}$ is an abelian
subalgebra $\z \subset {\gf}$ which is the centralizer of some
regular (but not necessarily semisimple) element  $g \in{\gf}$. In
particular, taking $g$ to be regular semisimple, we see that every
Cartan subalgebra (i.e. Lie algebra of a maximal torus)
is a regular centralizer. In fact, we will see in
\secref{regcent} that the set of regular
centralizers in ${\gf}$ is parametrized by an algebraic variety $\GNb$
which is a partial compactification of the parameter space $G/N$ for
the maximal tori. The simplest Higgs bundles are the {\it unramified}
ones, i.e. Higgs bundles $(E_G,\z_X)$ for which all the fibers of $\z_X$ are
maximal tori.

The situation is less transparent with spectral covers.
In fact, we do not know a good definition of a spectral
cover that would work for any $G$ and reproduce for $GL(n)$ the old
object.

Instead, we use the notion of a {\it cameral cover} introduced in
\cite{D1}. By definition, the latter is a finite flat map
$p:\tilX\to X$ such that the Weyl group $W$ of $G$ acts on $\tilX$ and
certain restrictions on the
ramification behaviour are satisfied (cf. \secref{intrHiggs}).
When $G=GL(n)$, we will note below that this notion is {\it different}
from that of a spectral cover, though equivalent to it.

It turns out that every Higgs bundle determines in a canonical way a
cameral cover, so one is led naturally to the problem of classification
of Higgs bundles with a given cameral cover. This is the problem we
solve in the present paper.
Given a cameral cover $\tilX$,
we will describe the corresponding Higgs bundles in terms of the
``Abelian'' data consisting of the maximal torus
$T\subset G$, the $W$-action on $T$, and the ramification pattern of
$\tilX$ over $X$. The ``non-Abelian'' data involving the group $G$
itself is not needed.

\subsection{Outline of the results} \label{0.2}

We formulate the above classification problem in
the categorical
framework, in terms of the
category $\Higgs_{\tilX}(X)$ of Higgs bundles together with
an isomorphism
between the induced cameral cover and $\tilX$.
Our first result shows that this classification problem is indeed abelian.

Namely, starting from $\tilX$ we define a sheaf of abelian groups $T_{\tilX}$.
We assert in \thmref{first} that $\Higgs_{\tilX}(X)$ is a {\it gerbe}
over the Picard
category of $T_{\tilX}$-torsors.
(These notions are reviewed for the reader's convenience in
\secref{gerbes}.)
This result has two immediate consequences.

First, the set of isomorphism classes of objects in our category
$\Higgs_{\tilX}(X)$, i.e. the set of isomorphism classes of Higgs
bundles with the given cameral cover $\tilX$,
if non-empty, carries a simply-transitive action of the abelian group
$H^1(X,T_{\tilX})$ (Corollary \ref{its_a_torsor}), and is therefore
non-canonically isomorphic to it. It is thus a generalized Prym variety,
cf. \cite{D2}: depending
on the circumstances, this may appear as a Jacobian of a spectral curve,
or as an ordinary Prym, or as various types of Prym-Tyurin varieties
\cite{K}, and so on.

The second consequence allows us to determine when Higgs
bundles with the given cameral cover $\tilX$ actually exist. This
happens if and only if the gerbe is trivial: the cameral cover $\tilX$
determines an obstruction class in $H^2(X,T_{\tilX})$, and Higgs
bundles with the given $\tilX$ exist if and only if this class vanishes
(Corollary \ref{obstruction_to_having_a_higgs}).

In the above, the sheaf $T_{\tilX}$ is defined in terms of the slightly larger sheaf
$\Tb_{\tilX}$ (on $X$) of $W$-equivariant maps $\tilX\to T$, i.e.
$\Tb_{\tilX}(U):=\Mor_W(\tilU,T)$, where $\tilU$
is the induced cameral cover of $U$.
For each positive root $\alpha: T \to \GG_m$, let $s_\alpha$ be the
corresponding reflection acting on $\tilX$, and let 
$D^\alpha_X \subset \tilX$
be its fixed point scheme. Any section $t$ of
$\Tb_{\tilX}(U)$ determines a function 
$\alpha \circ t: \tilU \to \GG_m$
which goes to its own inverse under the reflection
$s_ \alpha$. In particular, its restriction to the ramification locus
$D^{\alpha}_X$ {\em equals} its inverse, so it equals $\pm 1$.
The subsheaf $T_{\tilX} \subset \Tb_{\tilX}$
is given by the positive choice:
$$T_{\tilX}(U) :=\{ t \in \Tb_{\tilX}(U) \ | \
(\alpha \circ t) |_{ D^{\alpha}_U} = +1 \  \mbox{for each root} \
\alpha \}. $$

Although \thmref{first} is quite useful, it is not a completely
satisfactory result by itself, as it does not describe which
$T_{\tilX}$-gerbe we get. Our main result,
\thmref{second}, gives a complete description
of the category $\Higgs_{\tilX}(X)$ as the gerbe parametrizing certain
``$R$-twisted, $N$-shifted $W$-equivariant $T$-bundles on $\tilX$''.
The ``twist'' here is along the ramification divisors, and the ``shift''
is by the extension class of the normalizer $N$.

Our description of this gerbe is based on an
explicit description of the underlying Picard category
$\Tors_{T_{\tilX}}$ which
appears in the statement of \thmref{first}. An object in this category,
i.e. a $T_{\tilX}$-torsor, consists of:

\begin{itemize}

\item
A (weakly $W$-equivariant) $T$-bundle $\L_0$ on $\tilX$,

\item
A group homomorphism $\gamma_0:N_0\to \Aut(\L_0,\tilX/X)$, commuting
with the projctions to $W$, and

\item
For every simple root $\alpha_i$, a trivialization
$$\beta_{i,0}:\alpha_i(\L_0)|_{D^{\alpha_i}_X}\simeq
{\O}_{D^{\alpha_i}_X}.$$
\end{itemize}

The data of $\gamma_0$ and $\beta_0$ must satisfy some compatibility
conditions, which are described in detail in \secref{proof_second}.
(Roughly, these say that the collection $\beta_0$ of
isomorphisms $\beta_{i,0}$ is
$W$-equivariant, and $\beta_0, \gamma_0$ are related by the compatibility
condition: $\gamma_{0\vert D^{\alpha}_X} = \check{\alpha} \circ \beta_0.$)
Morphisms in this category are $T$-bundle maps that are
compatible with the data of $\gamma_0$ and $\beta_0$.

Our notation here is as follows. An element of the group
$\Aut(\L_0,\tilX/X)$,
for a $T$-bundle $\L_0$ on $\tilX$, consists of an element $w \in W$
together with an isomorphism $w^*(\L_0) \to \L_0$. The bundle
$\L_0$ is weakly $W$-equivariant if $\Aut(\L_0,\tilX/X)$
surjects onto $W$, in which case $\Aut(\L_0,\tilX/X)$ is an extension
of $W$ by $\Mor(\tilX, T)$. Now the semidirect product $N_0$ of $T$ and $W$
induces one such extension, and $\gamma_0$ is supposed to
induce an isomorphism of this extension with $\Aut(\L_0,\tilX/X)$. We think
of the root $\alpha$ as a homomorphism $T \to \GG_m$, so
$\alpha(\L_0)$ is the line bundle associated to $\L_0$ via this
homomorphism. Similarly, the coroot $ \check{\alpha}$ is a homomorphism
$\GG_m \to T$.

In describing our gerbe, we replace each linear feature in the
description of $\Tors_{T_{\tilX}}$ by an affine variant. We start with
the equivariance: the
$T$-bundles $\L_0$ were weakly $W$-equivariant (which means that
$w^*(\L_0)$ was isomorphic to $\L_0$, for each $w \in W$),
and in fact strongly $W$-equivariant (which simply means
that $W$ itself, and hence also the semidirect product
$N_0$, acted on them).

Our variant of the weakly
$W$-equivariant $T$-bundles $\L_0$ involves
$T$-bundles $\L$ which are $R$-twisted weakly $W$-equivariant, meaning
that now $w^*(\L) \otimes {\R}^w_X$ is isomorphic to $\L$, for each
$w \in W$. Here ${\R}^w_X$ is a $T$-bundle
on $\tilX$ which encodes the
ramification pattern of $\tilX$ over $X$. In the simplest case, when
$\tilX$ is integral and $w$ is the reflection $s_{\alpha}$ corresponding
to a simple root $\alpha$, we have
${\R}^w_X= {\R}^{\alpha}_X= \check{\alpha}(R^{\alpha}_X),$
where $R^{\alpha}_X$ is the line bundle ${\O}_{\tilX}(D^{\alpha}_X)$.
The precise definitions are given in \secref{main}.

Next, we need a substitute for the strong equivariance.
We replace $\Aut(\L_0,\tilX/X)$ by the group $\Aut_R(\L,\tilX/X)$
of isomorphisms $w^*(\L) \otimes {\R}^w_X \to \L$, and
the semidirect product $N_0$ by the normalizer $N$,
so we demand that $\gamma$ should map $N$ to $\Aut_R(\L,\tilX/X)$.

Finally, $\beta_i$ needs to be twisted by the ramification, so it now sends
$\alpha_i(\L)|_{D^{\alpha_i}_X} \to \R^{\alpha_i}|_{D^{\alpha_i}_X}$.
One final complication is that $\beta_i$ now depends (linearly) on the
choice of a lift of $w_i$ to an element $n_i \in N$. (This choice of
a lift is
unnecessary in the linear version, since $W$ is a subgroup of $N_0$, so
the $w_i$'s have a canonical lift.)

We can now give an almost complete statement of our main result,
\thmref{second}. It says that a Higgs bundle with given
cameral cover $\tilX$ is equivalent to:

\begin{itemize}

\item
An $R$-twisted, weakly $W$-equivariant $T$-bundle $\L$ on $\tilX$,

\item
A group homomorphism $\gamma:N\to \Aut_R(\L,\tilX/X)$, and

\item
For every simple root $\alpha_i$ and lift $n_i \in N$ of
the reflection $s_i\in W$ into $N$, the data of an isomorphism
$$\beta_i(n_i):\alpha_i(\L)|_{D^{\alpha_i}_X}\simeq
\R^{\alpha_i}|_{D^{\alpha_i}_X}.$$
\end{itemize}

The data of $\gamma$ and $\beta$ must satisfy several compatibility
conditions, which are described in detail in \secref{main}. (Roughly,
these say that the collection $\beta$ of isomorphisms $\beta_i(n_i)$ is
$N$-equivariant, and $\beta, \gamma$ are related by the compatibility
condition: $\gamma_{\vert D^{\alpha}_X} = \check{\alpha} \circ \beta.$)
In fact, the category $\Higgs_{\tilX}(X)$ is equivalent to the category
$\Higgs'_{\tilX}(X)$ whose objects are the triples $(\L, \gamma, \beta)$
as above.
Morphisms in this category are again $T$-bundle maps that are
compatible with the data of $\gamma$ and $\beta$.

Note that the possible non-triviality of our gerbe can be attributed to
three separate causes: the twist along the ramification $R$; the shift
resulting from non-triviality of the extension class of $N$; or the
extra complication involved in choosing the $\beta_i$. In
subsection \ref{simpleversion} we give a simplified version of our theorem,
which avoids this last complication. It applies
in all cases except when the group $G$ has $SO(2n+1)$ as a direct
summand.

\subsection{Some examples and applications} \label{0.3}

\sssec{The unramified case}
The cameral cover $\tilX \to X$ is unramified if and only if the Higgs
bundle $(E_G,\z_X)$ is unramified, i.e. if and only if the bundle of
regular centralizers $\z_X$ is actually a bundle of Cartan subalgebras.
In this case the classification (given in \cite{D2}) is easy: specifying a Higgs
bundle $(E_G,\z_X)$ with the unramified cameral cover $\tilX$ is
equivalent to giving an $N$-bundle $E_N$ over $X$ together with an
identification of the quotient $E_N/T$ with $\tilX$. In this case, our
$T$-bundle $\L$ is just $E_N$, considered as a $T$-bundle over $E_N/T =
\tilX$. Since there is no ramification, there is no
$R$-twist; similarly, there is no $\beta$; and $\Aut_R(\L,\tilX/X)$ is
just $\Aut(E_N,\tilX/X)$, which is
induced from the extension $N$, so $\gamma$ is the tautological map.

\sssec{$GL(n)$}

Consider first the case of $G=GL(n)$.
The spectral cover $\Xb$ is then of degree $n$ over $X$, while the cameral
cover $\tilX$ is of degree $n!$. The $n$
points of $\Xb$ above each point $x$ of $X$ correspond to the $n$
simultaneous eigenvectors (in the standard representation) of the
corresponding centralizer $\z_x$, while the $n!$ points of $\tilX$ above $x$
correspond to the ways of ordering these eigenvectors.
In a generic situation, e.g. when the Higgs bundle is unramified or
only simply ramified,
it is clear that $\tilX$ is precisely the Galois closure of the spectral
cover $\Xb$. Conversely, $\Xb$ is recovered as the quotient of $\tilX$
by $S_{n-1}$, the stabilizer in the permutation group
$W=S_{n}$ of one of the $n$ eigenvectors.
Following \cite{D1}, we study the relation between the two types of
covers in \secref{gln}. In particular, we show that the above
correspondence actually extends to an equivalence
between cameral and spectral covers, even when we are very far from the
generic situation.

\sssec{The universal objects}
The set of all maximal tori $T \subset G$, or equivalently the set of
Cartan subalgebras in $\gf$, is parametrized by the quotient $G/N$. Over
$X=G/N$ we have the tautological, unramified Higgs bundle: the underlying
$G$-bundle is the trivial one, $X \times G$, and the regular centralizers
are the universal family of Cartans. The corresponding (unramified)
cameral cover in
this case is $G/T \to G/N$. Note that a point of $G/T$ is determined by
a Cartan together with a Borel containing it.

The cover $G/T \to G/N$ admits a natural partial compactification
$\GTb \to \GNb$. Here $\GNb$ paramatrizes {\it regular centralizers}
in the Lie algebra $\gf$, and
$\GTb$ is the ramified $W$-cover of $\GNb$ parametrizing pairs consisting of a
regular centralizer together with a Borel containing it, cf. \secref{regcent}
and \secref{structure}. The map $\GTb \to \GNb$ is the cameral cover of
the tautologocal Higgs bundle on $\GNb$: the underlying
$G$-bundle is still $\GNb \times G$, and the regular centralizers
form the universal group-scheme $\Z$
of centralizers over $\GNb$.
We refer to these as universal objects; every
Higgs bundle on $X$ is locally the pullback of the tautological one
via some map $X \to \GNb$, and every cameral cover of $X$ is locally
the pullback of $\GTb \to \GNb$ via the same map $X \to \GNb$.

Although our ultimate results are concerned with Higgs bundles on arbitrary
schemes, much of our work boils down to a group-theoretic analysis of
these universal objects $\GNb$ and $\GTb$.
For instance, we will see that the ramification divisors are indexed by
the positive roots $\alpha$ of $G$. In fact, one
of the key points of this paper is that the tautological
group-scheme $\Z$ can be
completely recovered by looking
at the ramification pattern of $\GTb$ over $\GNb$. In a strong sense,
this says that a regular centralizer can be recovered from the scheme
parametrizing those Borels which contain it. This is our \thmref{centralizers}.
We emphasize that it is the phenomenon described in \thmref{centralizers}
which is ``responsible'' for the abelianization.

\sssec{$SL(2), PGL(2)$}
We saw that in the general case, the final form of the
answer is quite involved. A main source of technical difficulties is
the possible presense in $G$ of non-primitive coroots (cf. \cite{Sc}).
>From the classification of reductive groups we know that this
can occur only when $G$ has $SO(2n+1)$ as a direct factor. So the simplest
case where this extra complication occurs is for $G=SO(3)=PGL(2)$. In an
attempt to illustrate the effect of these non-primitive coroots,
we will, in \secref{examples}, work out explicitly and contrast the
examples of $G=SL(2)$, for
which no $\beta$'s are necessary because all coroots are primitive,
versus $G=PGL(2)$, for which the roots are non-primitive.
For these groups, both the spectral
cover and the cameral cover are double covers of $X$, so the entire
analysis can be made much more concrete than for a general group.
In particular, there are very explicit descriptions of the universal
objects $\GTb, \GNb, G/T, G/N$, cf. subsection\ref{rank1universals}.

\sssec{$K$-valued Higgs bundles}

The point of our abstract notion of a Higgs bundle is that it provides
a uniform approach to the analysis of various more concrete objects.
In the literature, the most common notion of a Higgs bundle is that of a
{\it $K$-valued Higgs bundle} on $X$, where $K$ is a fixed line bundle on
$X$. By definition, this means a pair $(E_G, s)$,
where $E_G$ is a principal $G$-bundle on $X$ and $s$ is a section
of $\gf_{E_G} \otimes K$.
Starting with one of our ``abstract'' Higgs bundles $(E_G,\z_X)$, we get a
$K$-valued Higgs bundle by choosing a section of $\z_X \otimes K$.
Conversely, a $K$-valued Higgs bundle $(E_G,s)$ on $X$ determines a
unique ``abstract'' Higgs bundle on the open subset $X_0 \subset X$ where
$s$ is regular. We say that a $K$-valued Higgs bundle is {\it regular}
if $X_0=X$.

Our philosophy is to think of a regular $K$-valued Higgs bundle as
involving two separate pieces of data.  The first requires specifying the
basis of ``eigenvectors'' of the Higgs field, i.e. it amounts to
specifying the underlying abstract Higgs bundle. The other piece of the
data corresponds to the ``eigenvalues''; in our case this amounts to
specifying the section $s$ of $\z_x\otimes K$. Our point is that
this second part of the data is irrelevant for the abelianization process,
so we focus on the ``eigenvectors'' encoded in the abstract Higgs bundle.
One obvious advantage of this approach is that it allows the bundle
$K$ of ``values" to be replaced by various other objects, as we will see
below.

A little more generally, we can work with the concept of a {\it regularized}
$K$-valued Higgs bundle on $X$, which means a triple
$(E_G,\z_X,s)$, with $(E_G,\z_X)$ a Higgs bundle in our
abstract sense, and $s$ a (not necessarily regular!) section
of $\z_X\otimes K$.
The moduli space of regular $K$-valued Higgs bundles is open
in the moduli of all $K$-valued Higgs bundles (for $X$ projective), and
is also open inside the moduli space of regularized $K$-valued
Higgs bundles.  For a "general" Higgs bundle, we can expect the
complement of
$X_0$ to have codimension 3, so if $X$ is projective of dimension 1 or 2,
we expect the open subset of regular Higgs bundles to be nonempty.

In Section \ref{values}, we apply
our results to show that the algebraic stack $\bfHiggs(X,K)$ of regularized
$K$-valued Higgs bundles on $X$ fibers over the affine space ${\mathbf B}(X,K)$ which
parametrizes $K$-valued cameral covers, i.e. pairs $(\tilX,v)$ where
$v$ is a $W$-equivariant map $v:\tilX\to \tf\otimes K$ (of
schemes over $X$). The fibers can be identified
with the gerbe $Higgs_{\tilX}(X)$ which we studied in the abstract case.
In accordance with our general philosophy,
the fiber is independent of the bundle $K$ or the way
$\tilX$ maps to $K$: it depends only on the abstract cameral cover
$\tilX$.

In case $X$ is a smooth, projective curve and $K$ is its canonical
bundle, we thus recover a version of
Hitchin's integrable system \cite{H}. (There is of course a difference,
in that we work with regularized $K$-valued Higgs bundles while Hitchin
uses semistable $K$-valued Higgs bundles.)
As an application, our results can be used to establish a duality between
the fibers of
the Hitchin map for a group $G$ and those corresponding to its
Langlands dual group $\check G$.

\sssec{Bundles on elliptic fibrations}

Essentially no new phenomena are encountered if we allow our Higgs
bundle to take its  ``values'' in a vector bundle $K$. But we can go
further and try to take $K$ to be any abelian group scheme over $X$,
such as the relative Picard scheme of some (projective, integral) family
$f:Y \to X$. This leads us in Section \ref{elliptics} to define
a {\it regularized $G$-bundle on $Y$} to be the data $(\tilX,E_G,\z_X)$,
with $\tilX \to X$ a cameral cover of $X$, and
$(E_G,\z_X) \in \Higgs_{\tilY}(Y)$ a Higgs bundle
on $Y$  with cameral cover $\tilY :=f^*\tilX$. This notion is most
natural in case $f$ is an {\it elliptic} fibration, since then
we know what it
means for a bundle (on $Y$) to be regular above a point (of $X$).
Just as was the situation for $K$-valued Higgs bundles, ``most''
$G$-bundles on an elliptic curve are indeed regular, and a regular bundle
has a unique regularization.

In \thmref{bundles on elliptics} we apply our results about abstract
Higgs bundles to obtain
a complete spectral description of regularized $G$-bundles on $Y$.
In the most interesting case, when $f$ is an elliptic fibration, this is
the main result of \cite{D3}. Letting $\bfReg(X,Y)$ denote the algebraic
stack of regularized $G$-bundles on $Y$, we obtain a ``spectral map''
$h: \bfReg(X,Y) \to {\mathbf B}(X,Y)$, sending a regularized bundle to its
$\Pic(Y/X)$-valued cameral cover, the fibers now being a slightly
twisted version of our gerbe $\Higgs_{\tilX}(X)$.

\subsection{Some history}\label{0.4}

The idea of abelianization has its source in quantum field theory and has
been extensively exploited by both physicists and mathematicians.
This idea was originally applied not to our notion of an abstract
Higgs bundle, but rather to $K$-valued Higgs bundles.
These were considered by Hitchin \cite{H} in case $X$ is a
curve and $K$ its canonical bundle. Other line bundles, on $X=P^1$,
were considered by Adams, Harnad and Hurtubise \cite{AHH} and Beauville
\cite{B}. Several aspects of spectral covers of $P^1$ and their Prym-Tyurin
varieties were considered by Kanev in \cite{K}.
The abelianization of $K$-valued Higgs bundles on other curves was
considered by Beilinson and Kazhdan, Bottacin, Donagi and Markman, Faltings,
Markman, and Scognamillo \cite{BK,Bo,DM,Fa,M,Sc}, among others.
In the case that the base $X$ is a curve, these Higgs bundles are
related to representations of the fundamental group of a punctured
Riemann surface, as well as to integrable systems arising from loop
algebras. The notion of a cameral cover was introduced in \cite{D1},
where its relation to the various spectral covers was analyzed.

The main point of many of the works cited above is to show, in various
interesting special cases, that the fiber of the Hitchin map, i.e. the
family of Higgs bundles with given spectral (or cameral) cover, "is"
generically a Jacobian or a Prym variety, depending on the group. A
description of this fiber in the general setting was announced in
\cite{D2}. In particular, the generalized Prym was described there as a
certain quotient of $H^1(\Tb_{\tilX})$. (This could be off by a finite
isogeny: we have seen that the correct description involves $H^1(T_{\tilX})$.)
It was also noted there that the fiber is canonically identified not with
the generalized Prym variety itself, but with a certain torsor over
it. The class of this torsor was described there in terms of the
"twist"  arising from the ramification divisor and the "shift"
by the class of the normalizer $N$ in $H^2(W,T)$. The additional
complication which arises only for $\on{SO}(2n+1)$ was first noted in
\cite{Sc}. This is encoded in the present work in our $\beta$'s.

Higgs bundles on higher dimensional varieties $X$, valued in the
cotangent bundle $K:=T^*X$, were introduced by Simpson \cite{Si}.
Through work of Corlette and Simpson, their moduli spaces are related to
those of local systems on $X$.
The version where $K$ is replaced by an elliptic fibration was
developed in \cite{D3} and \cite{Tanig}. These elliptically
valued Higgs bundles are of interest
because of their relevance to the construction and parametrization of
bundles on elliptic fibrations. These have attracted attention recently
because of their importance to understanding the conjectured duality
between F-theory and the heterotic string, cf.
\cite{FMW, FMW2, FMW3,FM,D3,CD,ICMP}.

\ssec{Notation} \label{0.5}
We work throughout with a fixed connected reductive group $G$
over $\CC$ and we let $\gf$ denote its Lie algebra.
We fix a Borel subgroup $B\subset G$ and denote by $\B$
the flag variety $G/B$. By definition, $\B$ classifies Borel
subalgebras in $\gf$.

Let $U$ be the unipotent radical of $B$ and $T$ the Cartan quotient $B/U$;
we will fix a splitting $T\to B$. We will denote by $\bof$ and
$\tf$ the Lie algebras of $B$ and $T$, respectively. The rank $r$ of
$G$ is by definition the dimension of $T$.
By $N$ we will denote the normalizer of $T$ (not the
nilpotent subgroup!), and by $W$ the Weyl group
$N/T$.

The set of positive roots will be denoted by $\Delta^+$. For
$\alpha\in \Delta^+$, let $\tf^\alpha\subset \tf$ denote the
corresponding root hyperplane and $s_\alpha\in W$ the corresponding reflection.
The set of simple roots we will denote by $I$. For $i\in I$, we will
use the notation $s_i$ instead of $s_{\alpha_i}$.

\part{Main results on Higgs bundles and cameral covers}

\section{Regular centralizers}  \label{regcent}

\ssec{} Recall that an element $x\in \gf$ is called {\it regular}
if its centralizer $Z_{\gf}(x)$ has the smallest
possible dimension, namely $r$ (the rank of $\gf$).
Note that with this definition, a regular element need not be semisimple.
The set of all regular elements forms an
open subvariety of $\gf$, which we will denote by $\gf_{reg}$.

A Lie subalgebra $\af\subset\gf$ is called a {\it regular centralizer} if
$\af=Z_{\gf}(x)$ for some $x\in \gf_{reg}$. Note that such $\af$
is automatically abelian.
Our first goal is to introduce a variety which parametrizes all regular
centralizers
in $\gf$.

\ssec{} Let $\on{Ab}^r$ be the closed sub-variety in the Grassmannian
of $r$-planes $Gr^r_{\gf}$
that classifies abelian subalgebras in $\gf$ of dimension
$r$. Let $\Gamma\subset \on{Ab}^r\times\gf$ be the incidence
correspondence, i.e. the closed subvariety defined by the condition:
$$(\af,x)\in\Gamma \text{ if } x\in\af.$$

Let $\Gamma_{reg}$ be the intersection
$\Gamma\cap(\on{Ab}^r\times\gf_{reg})$.

\begin{prop}  \label{GmodN}
There is a smooth morphism $\phi:\gf_{reg}\to \on{Ab}^r$ whose graph is
$\Gamma_{reg}$.
\end{prop}
The proof is postponed until \secref{structure}

\smallskip

Let $\GNb$ denote the image of the map $\phi$. The above proposition
implies that $\GNb$ is smooth and irreducible. It is clear that
$\CC$-points of $\GNb$ are exactly the regular centralizers in $\gf$.

\smallskip

By definition, the group $G$ acts on both $\on{Ab}^r$ and $\gf_{reg}$.
Therefore, the variety $\GNb$ acquires a natural $G$-action
and the map $\phi$ is $G$-equivariant.

\smallskip

Consider the quotient $G/N$; it classifies Cartan subalgebras in
$\gf$. These are the centralizers in $\gf$ of regular semisimple elements.
Hence, $G/N$
embeds into $\GNb$ as an open subvariety. Obviously, the action of $G$ on
$G/N$ by left multiplication is the restriction of its action on $\GNb$.

\ssec{} Consider the closed subvariety of $\GNb\times \B$ defined by
the condition: for $\af\in \GNb$ and $\bof'\in\B$,
$$(\af,\bof')\in \GTb \text{  if  }\af\subset \bof'.$$
We will denote this variety by $\GTb$ and the natural projection $\GTb\to
\GNb$
by $\pi$. It follows from the definitions that we have a natural
$G$-action on $\GTb$.

The quotient $G/T$ can clearly be identified with the open
sub-scheme $\pi^{-1}(G/N)$ of $\GTb$.
We have a natural action of the Weyl group $W=T\backslash N$ on $G/T$;
this action is
free and the quotient can be identified with $G/N$.

In what follows, by a $W$-cover of a scheme $X$
we will mean a finite flat scheme $p: \tilX\to X$, acted
on by $W$ such that $p_* {\O}_{\tilX}$
is locally isomorphic, as a coherent sheaf with a $W$-action,
to ${\O}_{X}\otimes \CC[W]$. A basic example is $\tof\to \tof/W$:
as is well known, it is ramified along the complexified walls of the
Cartan subalgebra $\tof$.

The following assertion will be proven in \secref{structure}

\begin{prop} \label{WactsonGmodT}
The variety $\GTb$ is smooth and connected. The $W$-action on $G/T$
extends
to the whole of $\GTb$ and it makes the latter a $W$-cover of
$\GNb$. Moreover, the two $W$-covers $\GTb\to \GNb$ and $\tof\to\tof/W$
are \'etale-locally isomorphic.
\end{prop}

\ssec{} \label{rank1universals}Here is an explicit description of $\GNb$ and $\GTb$ for
$G=SL(2)$.
In this case $\GNb$ is the space of all lines in $\gf$, i.e. $\GNb\simeq
\PP^2$.
We have a natural map $\GTb\to \PP^1\times\PP^1$, where the first
projection is
the natural map $\GTb\to\B\simeq \PP^1$ and the second projection is a
composition
of the first one with the action of $-1\in S_2\simeq W$ on $\GTb$.

It is easy to see that this map is an isomorphism. Under the
identification,
$\pi:\GTb\to\GNb$ is the symmetrization map $\PP^1\times\PP^1 \to \PP^2$.

\ssec{$G$-orbits}   \label{secorbits}
For each root $\alpha$, let $D^\alpha\subset\GTb$ denote the fixed point
set of $s_\alpha$ on $\GTb$. This is a smooth codimension $1$
subscheme of $\GTb$. Indeed, using the \'etale-local isomorphism
between $\GTb\to\GNb$ and $\tof\to\tof/W$ given in
\propref{WactsonGmodT}, it is enough to prove
this statement on $\tof$. However, $\tof^{s_\alpha}$
is just the corresponding root hyperplane $\tof^\alpha\subset\tof$.

\begin{prop} \label{orbits}
The $G$-orbits in $\GTb$ are precisely the locally closed subsets
$$D^{\Delta'}:=\underset{\alpha\in\Delta'}\cap (D^\alpha)\setminus
\underset{\beta\notin\Delta'}\cup (D^\beta)$$
where $\Delta'\subset \Delta$ is a subset of the set of roots, closed
under
linear combinations. The $G$-orbits in $\GNb$ are the images of the
$D^{\Delta'}$; they are indexed by the $\Delta'$ modulo the action of
$W$.
\end{prop}

The proof will be given in \secref{structure}.

\section{Higgs bundles and cameral covers}  \label{intrHiggs}

\ssec{Higgs bundles}
A family of Cartan subalgebras parametrized by a scheme $X$ is given by a
map from $X$
to $G/N$. Equivalently, it is given by a $G$-equivariant map from the
trivial $G$-bundle over $X$ to $G/N$. An advantage of this latter
description is that there is a natural way to twist it: given any
principal
$G$-bundle $E_G$ over $X$, we specify a
family of Cartan subalgebras {\it in the adjoint bundle}
$\gf_{E_G}:=E_G{\underset{G}\times}\gf$
by a $G$-equivariant map from $E_G$ to the variety $G/N$.
By generalizing this, we define:

\begin{defn} \label{DefHiggs}
A Higgs bundle over a scheme $X$ is a pair $(E_G,\sigma)$, where $E_G$ is
a
principal $G$-bundle over $X$ and $\sigma$ is a $G$-equivariant map
$\sigma:E_G\to \GNb$.
\end{defn}

\smallskip

Therefore, according to \propref{GmodN},
a Higgs structure in a given $G$-bundle $E_G$ is the same
as a vector subbundle $\z_X$ of $\gf_{E_G}$ of rank $r$ such that
$[\z_X,\z_X]=0$ and such that
locally in the \'etale topology $\z_X$ is the sheaf of centralizers of a
section of $E_G{\underset{G}\times}\gf_{reg}$.

\smallskip

The restriction of a Higgs bundle to an open subset $U \subset X$ over
which $E_G$ is trivialized can be specified more simply by a map $U \to
\GNb$.
In particular, the {\it universal Higgs bundle}
over $\GNb$ corresponds to the identity map $\GNb \to \GNb$.

\smallskip

\ssec{The Higgs category and stack}
Higgs bundles over $X$ form a category, denoted $\on{Higgs}(X)$.
An element of $\Hom((E^1_G,\sigma^1),(E^2_G,\sigma^2))$ is by
definition a $G$-bundle map $s:E^1_G\to E^2_G$ such that
$\sigma^2\circ s=\sigma^1$.

One can say that $\on{Higgs}(X)$ is the category of maps from $X$ to {\it the
stack} $G\backslash(\GNb)$.
Additionally, for a fixed $X$, we can consider the functor on the category of schemes, which
attaches to a scheme $S$ the category $\on{Higgs}(S\times X)$. When $X$ is projective,
this functor is representable by an algebraic stack, which we will denote by
${\bfHiggs}(X)$. (The representability follows because the stack ${\bfBun}_G(X)$
classifying principal $G$-bundles on $X$ is an algebraic stack. We have:
$G\backslash(\GNb)={\bfHiggs}(\on{Spec}(\CC))$.

\ssec{Cameral covers} \label{intr covers}
We will now introduce our second basic object.

\begin{defn} \label{Covers}
A $W$-cover of a scheme $X$ is a scheme $\tilX\overset{\pi}\to X$ finite and flat over $X$
such that as an $\O_X$-module with a $W$-action, $\pi_*(\O_{\tilX})$ is locally isomorphic to
$\O_X\otimes \CC[W]$.
\end{defn}

\begin{defn} \label{DefCam}
A cameral cover of $X$ is a $W$-cover $\tilX\to X$, such that locally with
respect to the \'etale topology on $X$, $\tilX$ is a pull-back of
the $W$-cover $\tof\to\tof/W$.
\end{defn}

As an example, we note that any $W$-cover is cameral when $G=\on{SL(2)}$,
i.e. $W=S_2$. On the other hand, not every $W=S_3$-cover is cameral: the
stabilizer of each point must be a Weyl subgroup of $W$, so, for example,
an $A_3$ stabilizer is not allowed.

\ssec{Openness} 
It is easy to see that the condition for a $W$-cover $\tilX\to X$
to be cameral is open on $X$.
Indeed, $\pi:\tilX\to X$ is cameral if and only if, locally on $X$, we can
find a $W$-equivariant embedding $\tilX\hookrightarrow X\times \tof$.
(Note that the space of $W$-equivariant maps of $X$-schemes $\tilX\to X\times \tof$
is ismorphic to the space of sections of the sheaf
$\Hom^W_{\O_X}(\tof^*\otimes \O_X,\pi_*(\O_{\tilX}))$, and the latter sheaf is
non-canonically isomorphic to $\tof\otimes \O_X$, since $\tilX\to X$ was
assumed to be a $W$-cover.)

\ssec{The cameral category and stack}
Cameral covers form a category in a natural way, denoted $\on{Cam}(X)$.
By definition, $\Hom(\tilX^1,\tilX^2)$ consists of all $W$-equivariant
isomorphisms $\tilX^1\to\tilX^2$.
It is easy to see that there exists an algebraic stack $\bfCam$, such that
$\on{Cam}(X)$ is the category $\Hom(X,\bfCam)$.

Indeed, consider the space of commutative $W$-equivariant
ring structures on the vector space $V:=\CC[W]$. This is clearly
an affine scheme, and let us denote it by $\bfCov$.
By construction, there exists a universal $W$-cover
$\widetilde{\bfCov}\to \bfCov$. Let $\bfCam'$ be the maximal
open subscheme of $\bfCov$, over which $\widetilde{\bfCov}$
is cameral. Let $Aut^W(V)$ be the algebraic group of automorphisms
of $V$ as a $W$-module. Clearly, $Aut^W(V)$ acts on $\bfCam'$
and the action lifts on $\widetilde{\bfCov}|_{\bfCam'}$.
We can now let $\bfCam$ be the stack-theoretic quotient
$Aut^W(V)\backslash \bfCam'$.

\medskip

As for Higgs bundles, for a fixed $X$ we can consider the functor
$S\mapsto \Cam(S\times X)$. For $X$ projective this functor is
representable by an algebraic stack $\bfCam(X)$.

\begin{prop}  \label{HiggsToCam}
There is a natural functor $F:\on{Higgs}(X)\to \on{Cam}(X)$.
In particular, for a projective scheme $X$, we obtain a map
between algebraic stacks $\bfHiggs(X)\to \bfCam(X)$.
\end{prop}

\begin{proof}

Any map $\sigma:  E_G \to \GNb$ determines a cameral cover
$\tilE_G$ of $E_G$, namely $\GTb{\underset{\GNb}\times}E_G$,
cf. \propref{WactsonGmodT}.

For a Higgs bundle, which involves a $G$-equivariant map $\sigma$,
the cameral cover $\tilE_G \to E_G$ is itself
$G$-equivariant, so by descent theory, it is pulled back from
a unique cameral cover $\tilX \to X$.

Clearly, the assignment $(E_G,\sigma)\mapsto \tilX$
constructed above is functorial.

\end{proof}

\medskip

Over an open set $U \subset X$ where $E_G$ is trivialized, the restriction
$\tilU \to U$ of the cameral cover is given in terms of $\sigma$ as:
$\GTb \times_{\GNb} U$. For example, applying this to the
universal Higgs bundle over $\GNb$ gives the cameral cover $\GTb
\to \GNb$. For this reason we refer in this paper to $\GTb \to
\GNb$ (rather than $\tof \to \tof/W$) as the {\it universal cameral
cover}.

\ssec{The fiber} Let us now fix a cameral cover $\tilX$. Let
$\on{Higgs}_{\tilX}(X)$
denote the category-fiber of the above functor
$F:\on{Higgs}(X)\to \on{Cam}(X)$ over $\tilX$.

In other words, the objects of $\on{Higgs}_{\tilX}(X)$ are pairs
$$((E_G,\sigma)\in\on{Higgs}(X),t:F(E_G,\sigma)\simeq\tilX)$$
and $\Hom((E^1_G,\sigma^1,t^1),(E^2_G,\sigma^2,t^2))$ is the set
of all bundle maps
$s:E^1_G\to E^2_G$ with $\sigma^2\circ s=\sigma^1$ and such that the
composition
$$\tilX{\overset{(t^1)^{-1}}\longrightarrow}F(E^1_G,\sigma^1)\to
F(E^2_G,\sigma^2)  {\overset{t^2}\longrightarrow}\tilX$$ is the identity
automorphism of $\tilX$.

The goal of this paper
is to describe explicitly the category $\on{Higgs}_{\tilX}(X)$
in terms of the $W$-action on $\tilX$.

\section{Gerbes}  \label{gerbes}

\ssec{}
Since the objects we study have automorphisms, it is difficult to
describe them adequately without the use of some categorical language.
Specifically, our description requires the notion of an $A$-gerbe, where
$A$ is a sheaf of abelian groups on $X$.  This is a
particularly useful case of the more general notion of a gerbe over a
sheaf of Picard categories. In this section we review
the corresponding definitions. For more details, the reader is
referred to \cite{Gi} or \cite{DeMi}

Let $\et(X)$ denote the big \'etale site over $X$. (By definition,
$\et(X)$
is the category of all schemes over $X$ and the covering maps are
surjective
\'etale morphisms.)

\ssec{} Recall that a presheaf $\C$ of categories on $\et(X)$ assigns
to every object $U\to X$ in $\et(X)$ a category $\C(U)$ and to every
morphism $f:U_1\to U_2$ in $\et(X)$ a functor $f^*_\C:\C(U_2)\to\C(U_1)$.
Moreover, for every composition $U_1{\overset{f}\to} U_2{\overset{g}\to}
U_3$ there should be a natural transformation $f^*_\C\circ
g^*_\C\twoheadrightarrow (g\circ f)^*_\C$, such that an obvious
compatibility relation for three-fold compositions holds.

\medskip

A presheaf $\C$ of categories on $\et(X)$ is said to be
a sheaf of categories (or a stack) if the following two axioms hold:

\smallskip

\noindent {\bf Axiom SC-1.} For $U\to X$ in $\et(X)$ and a pair of objects
$C_1,C_2\in \C(U)$, the presheaf of sets on $\et(U)$
that assigns to $f:U'\to U$ the set
$$\Hom_{\C(U')}(f_\C^*(C_1),f_\C^*(C_2))$$ is a sheaf.

\smallskip

\noindent {\bf Axiom SC-2.} If $f:U'\to U$ is a covering,
then the category $\C(U)$ is equivalent to
the category of descent data on $\C(U')$ with respect to $f$ (i.e. every
descent data on $\C(U')$ with respect to $f$ is canonically effective,
cf. \cite{DeMi}, p. 221).

\ssec{} Here is our main example of a sheaf of categories. Fix a cameral
cover
$\tilX\to X$. For every object $U\in\et(X)$ write
$\tilU:=U\underset{X}\times \tilX$,
which is a cameral cover of $U$.

We define the presheaf of categories $\Higgs_{\tilX}$ by
$\Higgs_{\tilX}(U):=\Higgs_{\tilU}(U)$ (the functors $\Higgs_{\tilX}(U)\to
\Higgs_{\tilX}(U')$
for $U'\to U$ and the corresponding natural transformations are defined in
a natural way).

The following is an easy exercise in descent theory:

\begin{lem}
$\Higgs_{\tilX}$ satisfies SC-1 and SC-2.
\end{lem}

\ssec{} Recall that a Picard category is a groupoid endowed with a a
structure of a
tensor category, in which every object is invertible. A basic example (and
the
source of the name) is the category of line bundles over a scheme.

\smallskip

A sheaf of categories $\P$ is said to be a sheaf of Picard categories if
for every $(U\to X)\in\et(X)$, $\P(U)$ is endowed with a structure of a
Picard category such that the pull-back functors $f^*_\P$ are
compatible with the tensor structure in an appropriate sense. If $\P_1$
and $\P_2$ are two sheaves of Picard categories one defines (in a
straightforward fashion) a notion of a tensor functor between them.

\smallskip

A typical and the most important example of a sheaf of Picard categories
can be constructed as follows:

Let $\A$ be a sheaf of abelian groups over $\et(X)$. For an object $f:U\to
X$ of
$\et(X)$ let $\Tors_\A(U)$ denote the category
of $\A|_U$-torsors on $U$. This is a Picard category and it is easy to
see that the assignment $U\to \Tors_\A(U)$ defines a sheaf of Picard
categories on $\et(X)$ which we will denote by $\Tors_\A$.

\ssec{} Just as a torsor is a space on which a group acts simply
transitively, a gerbe is a category on which a Picard category acts simply
transitively:
A category $\C$ is said to be a {\it gerbe over the Picard category}
$\P$, if $\P$ acts on $\C$ as a tensor category and for any object $C\in
\C$ the functor $\P\to\C$ given by $$P\in\P\Longrightarrow
\on{Action}(P,C)\in\C$$ is an equivalence.

\smallskip

Now, if $\P$ is a sheaf of Picard categories and $\C$ is another
sheaf of categories we say that $\C$ is a {\it gerbe over the sheaf
of Picard categories} $\P$, if the following holds:

\begin{itemize}

\item
For every $(U\to X)\in\et(X)$, $\C(U)$ has
a structure of a gerbe over $\P(U)$. This structure is compatible with
the
pull-back functors $f^*_\P$ and $f^*_\C$.

\item
There exists a covering $U\to X$, such that $\C(U)$ is non-empty.

\end{itemize}

A basic feature of gerbes is that if $\C_1$ and $\C_2$ are gerbes over 
$\P$, one can form a new gerbe $\C_1\underset{\P}\otimes \C_2$,
called their tensor product, cf. \cite{DeMi}.

\ssec{}\label{Q}
The basic example of a gerbe over an arbitrary sheaf of Picard
categories $\P$, is $\P$ itself. Here is a less trivial example:

Fix a short exact sequence $0 \to \A \to \A'' \to \A' \to 0$ of sheaves of
abelian groups on $X$ and let $\tau_{\A'}$ be an $\A'$-torsor over $X$.
We introduce a sheaf of categories $\C=\C_{\tau_{\A'}}$ as follows.
For $U\in\et(X)$, $\C(U)$ is the category of
all ``liftings'' of $\tau_{\A'}|_U$ to an $\A''|_U$-torsor. It is easy to
check that
$\C$ is a gerbe over $\P=\Tors_\A$.

In fact, gerbes over $\Tors_\A$ can be classified
cohomologically:

\begin{lem}  \label{cohomologyclass}
There is a bijection between the set of equivalence classes of gerbes
over $\Tors_\A$ and $H^2(X,\A)$.  For a given gerbe $\C$ the
corresponding class in $H^2(X,\A)$ vanishes if and only if the category
$\C(X)$ of "global sections" is non-empty.
\end{lem}

In the above example, the class on $H^2(X,\A)$ corresponds to the image of
the class
of $\tau_{\A'}$ under the boundary map $H^1(X,\A') \to H^2(X,\A)$.

\ssec{}  \label{exact sequences}
The following will be needed in \secref{elliptics}.

Let ${\mathbf a}:\P_1$ and $\P_2$ be sheaves of Picard categories, and
$\P_1\to \P_2$ be a functor compatible with the tensor structure.
We say that ${\mathbf a}$ is a monomorphism if for every $U\in \et(X)$
the functor ${\mathbf a}(U):\P_1(U)\to \P_2(U)$ is faithful.

We say that ${\mathbf a}$ is an epimorphism if
for every $U\in \et(X)$ and $P,\widetilde{P}\in \P_1(U)$ the map of sheaves
on $\et(U)$:
$\Hom_{\P_1(U')}(P|_{U'},\widetilde{P}|_{U'})\to
\Hom_{\P_2(U')}({\mathbf a}(P|_{U'}),{\mathbf a}(\widetilde{P})|_{U'})$
is an epimorphism (in the sense of sheaves).

Similary, if we have three sheaves of Picard categories and tensor
functors ${\mathbf a}:\P_1\to \P_2$, and ${\mathbf b}:\P_2\to \P_3$
we say that the form a short exact sequence if ${\mathbf b}$ is an epimorphism
and ${\mathbf a}$ induces an equivalence between $\P_1$ and the category-fiber
of $\P_2$ over the unit object in $\P_3$.

In this case, for every object $P_3\in \P_3$, the category fiber of
$\P_2$ over it is, in a natural way, gerbe over $\P_1$. This generalizes
the above example of $0 \to \A \to \A'' \to \A' \to 0$.

\medskip

Let now $\C_1$ be a gerbe over $\P_1$, and ${\mathbf a}:\P_1\to \P_2$ a tensor functor.
In this case one can construct a canonical induced gerbe $\C_2$ over $\P_2$
with the property that there exists a functor $\C_1\to \C_2$,
compatible with the $\P_1$- and $\P_2$-actions via ${\mathbf a}$.

Suppose now that
$$0\to \P_1\overset{{\mathbf a}}\to \P_2\overset{{\mathbf b}}\to \P_3\to 0$$ is a short exact
sequence of Picard categories, and $\C_1$ is a gerbe over $\P_1$.
Let $\C_2$ be the corresponding induced $\P_2$-gerbe.

The next lemma follows from the definitions in a straightforward way:

\begin{lem} \label{induced gerbes}
There exists a canonical functor $\C_2\to\P_3$. The category fiber of
$\C_2$ over a given object $P_3\in \P_3$ is naturally a $\P_1$-gerbe,
canonically equivalent to the tensor product $\C_1\underset{\P_1}\otimes {\mathbf b}^{-1}(P_3)$.
\end{lem}

\section{$\Higgs_{\tilX}$ is a gerbe} \label{HiggsIsGerbe}

\ssec{}

Given a cameral cover $\tilX\to X$, let $\Tb_{\tilX}$ be the sheaf
``of $W$-equivariant maps $\tilX\to T$'' on the \'etale site over $X$.
More precisely,
for $U\in\et(X)$, $\Tb_{\tilX}(U)=\Hom_W(\tilU,T)$, where $\tilU$
is the induced cameral cover of $U$
and the subscript ``$W$'' means maps respecting the $W$-action.

However, we need a slightly smaller sheaf.

\ssec{} \label{sheaf}

Let $D^\alpha_X$ (for each positive root $\alpha$) be the fixed point
scheme of the
reflection $s_\alpha$ acting on $\tilX$. Locally, this is
the pullback of the {\it universal ramification divisor}, i.e.
$D^\alpha \subset \GTb$.

\smallskip

Let $\alpha$ be a root of $G$, considered as a homomorphism
$\alpha: T \to \GG_m$. Then any section $t$ of
$\Tb_{\tilX}(U)$ determines a function $\alpha \circ t:U\to \GG_m$
which goes to its own inverse under the reflection
$s_ \alpha$. In particular, its restriction to the ramification locus
$D^{\alpha}_X$ {\em equals} its inverse, so it equals $\pm 1$.
The subsheaf $T_{\tilX} \subset \Tb_{\tilX}$
is defined by the following condition:
$$T_{\tilX}(U) :=\{ t \in \Tb_{\tilX}(U) \ | \
(\alpha \circ t) |_{ D^{\alpha}_U} = +1 \  \mbox{for each root} \
\alpha \}. \;\;\;\;\;\;\;\;\; (*)$$

By construction, $\Tb_{\tilX}/T_{\tilX}$ is a
${\mathbb Z}_2$-torsion sheaf. Note, in addition, that it suffices to
impose
condition (*) for one representative of each orbit of
$W$ on the set of roots.

\smallskip

\noindent{\bf Remark.}
Recall that a coroot $\check\alpha:\GG_m\to T$ is called primitive if
$\on{ker}(\check\alpha)=1$
(this is equivalent to saying that $\check\alpha$ is a primitive element
of the lattice of cocharacters of $T$.) It is clear that condition (*)
holds automatically
for roots whose corresponding coroots are primitive. For example, when the
derived group of $G$ is simply-connected, all coroots
are primitive, i.e. (*) is automatic and $T_{\tilX}=\Tb_{\tilX}$. In fact,
$G$ has
non-primitive coroots if and only if it contains $SO(2n+1)$ (e.g.
$PGL(2)=SO(3)$)
as a direct factor, as is easily seen from the classification of Dynkin
diagrams.

\ssec{} Our first result can be stated as:

\begin{thm} \label{first}
$\Higgs_{\tilX}$ is a gerbe over $\Tors_{T_{\tilX}}$.
\end{thm}

Let us list several corollaries of this theorem:

\begin{cor}\label{obstruction_to_having_a_higgs}
To a cameral cover $\tilX$ there corresponds a class in
$H^2(X,T_{\tilX})$,
which vanishes if and only if $\tilX$ is the cameral cover corresponding
to some Higgs bundle.
\end{cor}

This is immediate from \lemref{cohomologyclass}.

\begin{cor}\label{its_a_torsor}
Suppose $\Higgs_{\tilX}(X)$ is non-empty. The set of isomorphism classes
of
objects in this category carries a simply-transitive action of
$H^1(X,T_{\tilX})$.
The group of automorphisms of every object is canonically isomorphic to
$T_{\tilX}(X)$.
\end{cor}

\section{Ramification}  \label{secramification}

\ssec{} We now proceed to the formulation of our main result,
\thmref{second},
which describes the category
$\Higgs_{\tilX}(X)$ completely in terms of $\tilX$.
For that purpose, we need to introduce some further
notation that has to do with the ramification pattern of $\tilX$ over $X$.

\ssec{} For each root $\alpha$ we will define a line bundle $R_X^\alpha$
on $\tilX$. Assume first that $\tilX$ is integral. In this case
the subscheme $D_X^{\alpha}\subset \tilX$ is a Cartier divisor, because
locally it is the pull-back of $D^\alpha\subset\GTb$. We set
$R^\alpha_X=\O(D_X^\alpha)$.

When $\tilX$ is arbitrary we proceed as follows. The construction is
local,
so we may assume that $X$, and hence also $\tilX$, is affine.
Let $I^\alpha_X$ be a coherent
sheaf on $\tilX$ generated by symbols $\{g\}$, for
$ \{ g\in \O_{\tilX}\,|\, s_\alpha(g)=-g \} $ that satisfy the relations:
$$f\cdot \{g\}=\{f\cdot g\} \text{ for all } f \text{ such that }
s_\alpha(f)=f.$$

Locally, $I^\alpha_X$ is the pull-back of the sheaf of ideals of the
subscheme
$D^\alpha\subset \GTb$. Hence, $I^\alpha_X$ is a line bundle. We have a
natural map
$I^\alpha_X\to \O_{\tilX}$ that sends $\{g\}\mapsto g$ and, by
construction,
its cokernel is $\O_{D^\alpha_X}$.

We define the line bundle $R^\alpha_X$ as the inverse of $I^\alpha_X$.
We have a canonical section $\O_{\tilX}\to R^\alpha_X$ whose locus of
zeroes
is the subscheme $D^\alpha_X$.

\ssec{} Consider the $T$-bundle $\R^\alpha_X:=\check\alpha(R^\alpha_X)$
(i.e. $\R^\alpha_X$ is induced from $R^\alpha_X$ by means of
the homomorphism $\check\alpha:\GG_m\to T$).

\smallskip

For an element $w\in W$ we introduce the $T$-bundle $\R^w_X$ on $\tilX$ as
$$\R^w_X:={\underset{\alpha}\otimes}\, \R_X^{\alpha},$$
\noindent where $\alpha$ runs over those positive roots for which
$w(\alpha)$ is negative. For example, for $w=s_i$ (a simple reflection),
$\R^{s_i}_X\simeq \R^{\alpha_i}_X$.

\smallskip

Observe that given a $T$-bundle $\L$ on $\tilX$ and an element $w\in W$,
there are
two ways to produce a new $T$-bundle: we can pull back by $w$ acting as an
automorphism of $\tilX$, or we can conjugate the $T$-action by $w$. We
will {\em always} write $w^*(\L)$ for the combination of {\em both}
actions. For example, for $G=SL(2)$, the $T$-bundle $\L$ is equivalent
to a line bundle $L$. The two individual actions on $\L$ of the
non-trivial element $-1\in S_2=W$ send $L$ to $(-1)^*(L)$ and $L^{-1}$,
respectively, while $(-1)^*(\L)$ corresponds to the line bundle
$(-1)^*(L^{-1})$.
In particular, we have:
\begin{equation} \label{ramroots}
w^*(\R^\alpha_X)\simeq \R^{w^{-1}(\alpha)}_X.
\end{equation}

\begin{lem} \label{ramcocycle}
There is a canonical isomoprhism
$\R^{w_1\cdot w_2}_X\overset{\varpi(w_1,w_2)}\longrightarrow
w_2^*(\R^{w_1}_X)\otimes \R^{w_2}_X$.
\end{lem}

The proof follows imediately from the definition of $\R^w_X$ and \eqref{ramroots}.
The following proposition is necessary
for the formulation of \thmref{second}.

\begin{prop} \label{Rtriviality}
Let $\alpha_i$ be a simple root and let $w\in W$ be such that
$w(\alpha_i)=\alpha_j$ (another positive simple root). Then, the line
bundle
$\alpha_i(\R_X^w)|_{D^{\alpha_i}_X}$ admits a canonical trivialization.
\end{prop}

\begin{proof}

Let us observe first that, since we are using only roots rather than
arbitrary weights, it is sufficient to consider the case when $[G,G]$ is
simply-connected.

We have $w\cdot s_i=s_j\cdot w$, hence, by \lemref{ramcocycle}
$$s_i^*(\R^w_X)\otimes \R^{s_i}_X\simeq \R^{w\cdot s_i}_X\simeq
\R^{s_j\cdot w}_X\simeq w^*(\R^{s_j}_X)\otimes \R^w_X.$$

However, by definition $w^*(\R^{s_j}_X)\simeq \R^{s_i}_X$, so
we obtain that $s_i^*(\R^w_X)\simeq \R^w_X$. By restricting to
$D^{\alpha_i}_X$,
we obtain $\check\alpha_i(\alpha_i(\R^w_X))\simeq
\check\alpha_i(\O_{D^{\alpha_i}_X})$.

Since $[G,G]$ is simply-connected, every coroot is primitive. Therefore,
there
exists a weight $\lambda$, such that
$\lambda\circ\check\alpha_i=\on{id}:\GG_m\to \GG_m$.
By applying $\lambda$ to the above isomorphism
$\check\alpha_i(\alpha_i(\R^w_X))\simeq
\check\alpha_i(\O_{D^{\alpha_i}_X})$, we obtain
an isomorphism $\alpha_i(\R^w_X)\overset{\on{isom}^\lambda}\simeq
\O_{D^{\alpha_i}_X}$.

Now it only remains to check that this isomorphism is independent of the
choice of $\lambda$.
However, since the $\R^w_X$'s are locally pull-backs of the corresponding
$T$-bundles on $\GTb$,
it suffices to consider the universal situation, namely the case $X=\GNb$.

In the latter case, the $T$-bundle $\R^w|_{D^{\alpha_i}}$ itself is
trivialized over
an open dense part of $D^{\alpha_i}$, namely over $D^{\alpha_i}-
\underset{\alpha\neq\alpha_i}\cup (D^\alpha\cap D^{\alpha_i})$. This is
because $\alpha_i$
is not among the set of roots which become negative under the action of
$w$. In
particular, we obtain an isomorphism
$\alpha_i(\R^w)\overset{\on{isom}'}\simeq \O_{D^{\alpha_i}}$ over
$D^{\alpha_i}-\underset{\alpha\neq\alpha_i}\cup (D^\alpha\cap
D^{\alpha_i})$.

Moreover, it is easy to see that for any $\lambda$ as above, the
isomorphisms
$\on{isom}^\lambda$ and $\on{isom}'$ coincide. In particular,
$\on{isom}^\lambda$ is independent of $\lambda$ over
$D^{\alpha_i}-\underset{\alpha\neq\alpha_i}\cup (D^\alpha\cap
D^{\alpha_i})$
and hence over the whole of $D^{\alpha_i}$, which is what we need.

$$  $$

\end{proof}

The following notions will be used in the formulation of \thmref{second}.

\begin{defn}
Let $\L_0$ be a $T$-bundle on $\tilX$. We say that it is weakly
$W$-equivariant if for every $w$ there exists an isomorphism
$w^*(\L_0)\to \L_0$.
\end{defn}

For a weakly $W$-equivariant $T$-bundle, let $\on{Aut}(\L_0)$ be the group
whose elements are pairs: an element $w\in W$ plus an isomorphism
$w^*(\L_0)\to \L_0$. By definition,
$\on{Aut}(\L_0)$ fits into a short exact sequence:
$$1\to \Hom(\tilX,T)\to \on{Aut}(\L_0)\to W\to 1.$$

\begin{defn}
A strongly $W$-equivariant $T$-bundle is a
weakly $W$-equivariant $T$-bundle $\L_0$ plus a choice
of a splitting $\gamma_0:W\to \on{Aut}(\L_0)$.
\end{defn}

\medskip

\begin{defn}
A $T$-bundle on $\tilX$ is called weakly $R$-twisted $W$-equivariant
if for every $w\in W$ there exists an isomorphism
$w^*(\L)\otimes \R^w_X\simeq \L$.
\end{defn}

For a weakly $R$-twisted $W$-equivariant $T$-bundle $\L$ we introduce the
group
$\on{Aut}_\R(\L)$. Its elements are pairs $w\in W$ and an ismorphism
$w^*(\L)\otimes \R^w_X\simeq \L$. The group law is defined via the
isomorphism
$\varpi(w_1,w_2)$ of \lemref{ramcocycle}.
By definition, $\on{Aut}_\R(\L)$
is also an extension of $W$ by means of $\Hom(\tilX,T).$

\section{The main result}\label{main}

\ssec{} We need one more piece of notation. For a simple root
$\alpha_i$, let $M_i$ be the corresponding minimal Levi subgroup.
Under the projection $N\to W$, the intersection $N\cap [M_i,M_i]$
surjects onto $\langle s_i\rangle\simeq S_2$. Let $\N_i$ denote
the preimage of $s_i$ in $N\cap [M_i,M_i]$.

By definition,
if $n_i$ and $n'_i$ are two elements in $\N_i$, there exists
$c\in \GG_m$ such that $n'_i=\check\alpha_i(c)\cdot n_i$.

\ssec{} Given a cameral cover $\tilX\to X$, we introduce the category
$\Higgs'_{\tilX}(X)$ of ``$R$-twisted, $N$-shifted $W$-equivariant
$T$-bundles on $\tilX$''. Its objects consist of:

\begin{itemize}

\smallskip

\item
A weakly $R$-twisted $W$-equivariant $T$-bundle $\L$ on $\tilX$.

\smallskip

\item
A map of short exact sequences:
$$
\CD
1 @>>> T @>>> N @>>> W @>>> 1 \\
@VVV  @V{\text{natural map}}VV     @V{\gamma}VV   @V{\on{id}}VV   @VVV \\
1 @>>> \on{Hom}(\tilX,T) @>>> \on{Aut}_\R(\L) @>>> W @>>> 1
\endCD
$$

\smallskip

\item

For each simple root $\alpha_i$ and element $n_i\in \N_i$,
an isomorphism of line bundles on $D^{\alpha_i}_X$
$$\beta_i(n_i):\alpha_i(\L)|_{D^{\alpha_i}_X}\simeq
R^{\alpha_i}_X|_{D^{\alpha_i}_X}.$$

\end{itemize}

These data must satisfy three compatibility conditions:

\smallskip

\noindent (1) If $n'_i=\check\alpha_i(c)\cdot n_i$ for $c\in {\mathbb
G}_m$,
then $\beta_i(n'_i)=c\cdot \beta_i(n_i)$.

\medskip

\noindent (2) Let $\alpha_i$ be again a simple root and $n_i\in \N_i$.
Consider the isomorphism
$$\gamma(n_i):s^*_i(\L)\otimes \R^{s_i}_X\simeq \L.$$
When we restrict it to $D^{\alpha_i}_X$ it induces an isomorphism
$$\check\alpha_i(\alpha_i(\L)|_{D^{\alpha_i}_X})\simeq
\check\alpha_i(R^{\alpha_i}_X|_{D^{\alpha_i}_X}),$$
by the definition of $\R^{s_i}_X$.
We need that this isomorphism coincides with
$\check\alpha_i(\beta_i(n_i))$.

\medskip

\noindent (3)  Let $\alpha_i$ and $\alpha_j$ be two simple roots and let
$w\in W$ be such that $w(\alpha_i)=\alpha_j$. Let $\tilde{w}\in N$
be an element that projects to $w$, and $n_j$ be an element of $\N_j$. By
pulling back
the isomorphism $\beta_j(n_j)$ with respect to $w$,
we obtain an isomorphism $\alpha_i(w^*(\L))|_{D^{\alpha_i}_X}\simeq
R^{\alpha_i}_X|_{D^{\alpha_i}_X}$.
In addition, the isomorphisms induced by $\gamma(\tilde{w})$ and
\propref{Rtriviality} lead
to a sequence of isomorphisms:
$$\alpha_i(\L)|_{D^{\alpha_i}_X}\overset{\gamma(\tilde{w})}\longrightarrow
\alpha_i(w^*(\L))|_{D^{\alpha_i}_X}\otimes
\alpha_i(\R^w_X)|_{D^{\alpha_i}_X}
\overset{\text{\propref{Rtriviality}}}\longrightarrow
\alpha_i(w^*(\L))|_{D^{\alpha_i}_X}.$$

By composing the two, we obtain an isomorphism
$\alpha_i(\L)|_{D^{\alpha_i}_X}\simeq R^{\alpha_i}_X|_{D^{\alpha_i}_X}$
and our condition is that it coincides with $\beta_i(n_i)$, where
$n_i=\tilde{w}^{-1}\cdot n_j\cdot \tilde{w}\in \N_i$.

\medskip

This concludes the definition of objects of $\Higgs'_{\tilX}(X)$.
Morphisms between $(\L,\gamma,\beta_i)$ and
$(\L^1,\gamma^1,\beta^1_i)$ are $T$-bundle isomorphism maps $\L^1\to\L$,
which intertwine in the obvious sense $\gamma$ with $\gamma^1$ and
$\beta_i$ with $\beta^1_i$.

\ssec{} It is easy to see that
$\Higgs'_{\tilX}(X)$ can be naturally sheafified. Namely, we define the
presheaf
of categories $\Higgs'_{\tilX}$ by setting for for $U\in\et(X)$,
$\Higgs'_{\tilX}(U):=\Higgs'_{\tilU}(U)$.
The pull-back functors are defined in an evident manner and it is easy to
see that $\Higgs'_{\tilX}$
satisfies SC-1 and SC-2.

Our main result is:

\begin{thm} \label{second}
The sheaves of categories $\Higgs'_{\tilX}$ and $\Higgs_{\tilX}$ are
naturally equivalent.
\end{thm}

In particular, we obtain that $\Higgs_{\tilX}(X)$ is equivalent to
$\Higgs'_{\tilX}(X)$. In other words, a Higgs bundle on $X$ with the
given cameral cover $\tilX$ is equivalent to a $T$-bundle on $\tilX$
which is $R$-twisted, $N$-shifted $W$-equivariant.

\ssec{Variant}  \label{simpleversion}

Assume that all coroots in $G$ are primitive, i.e. for every $\alpha$,
the corresponding $1$-parameter subgroup maps injectively into $T$.

We claim that the definition of $\Higgs'_{\tilX}(X)$ is equivalent
to the following (simplified) one. We introduce the category
$\Higgs''_{\tilX}(X)$
as follows:

Objects of $\Higgs''_{\tilX}(X)$ are pairs

\begin{itemize}

\smallskip

\item
A weakly $R$-twisted $W$-equivariant $T$-bundle $\L$ on $\tilX$.

\smallskip

\item
A map of short exact sequences:
$$
\CD
1 @>>> T @>>> N @>>> W @>>> 1 \\
@VVV  @V{\text{natural map}}VV     @V{\gamma}VV   @V{\on{id}}VV   @VVV \\
1 @>>> \on{Hom}(\tilX,T) @>>> \on{Aut}_\R(\L) @>>> W @>>> 1,
\endCD
$$

\end{itemize}
such that the following condition holds:

\smallskip

\noindent (1') Let $\lambda$ be a weight of $T$ such that
$\langle\lambda,\check\alpha_i\rangle=0$, which implies that
$\lambda(\L)|_{D^{\alpha_i}_X}\simeq\lambda(s_i^*(\L)\otimes
\R^{s_i}_X)|_{D^{\alpha_i}_X}$.
Our condition is that for every $n_i\in \N_i$ the composition
$$\lambda(\L)|_{D^{\alpha_i}_X}\simeq
\lambda(s_i^*(\L)\otimes
\R^{s_i}_X)|_{D^{\alpha_i}_X}\overset{\gamma(n_i)}
\longrightarrow \lambda(\L)|_{D^{\alpha_i}_X}$$
is the identity map.

Morphisms between $(\L,\gamma)$ and $(\L^1,\gamma^1)$ are $T$-bundle maps,
which intertwine between $\gamma$ and $\gamma^1$.

\medskip

Let us show that $\Higgs'_{\tilX}(X)$ and $\Higgs''_{\tilX}(X)$ are
naturally equivalent.
Indeed, if we have an object $(\L,\gamma,\beta_i)\in\Higgs'_{\tilX}(X)$,
the corresponding object of $\Higgs''_{\tilX}(X)$ is obtained by just
forgetting the $\beta_i$'s.

\smallskip

Conversely, if $(\L,\gamma)\in \Higgs''_{\tilX}(X)$, we reconstruct
the $\beta_i$'s as follows:

For a simple root $\alpha_i$ and $n_i\in\N_i$ consider the isomorphism
$\gamma(n_i)$ restricted to $D^{\alpha_i}_X$. It yields an isomorphism
$$\check\alpha_i(\alpha_i(\L))|_{D^{\alpha_i}_X}\simeq
\check\alpha_i(R^{\alpha_i}_X)|_{D^{\alpha_i}_X}.$$

Since $\check\alpha_i$ is primitive, there exists a weight
$\lambda'$ with $\langle\lambda',\check\alpha_i\rangle=1$. By evaluating
$\lambda$ on the above
isomorphism, we obtain the required identification
$\beta_i(n_i):\alpha_i(\L)|_{D^{\alpha_i}_X}\to
R^{\alpha_i}_X|_{D^{\alpha_i}_X}$.
This isomorphism does not depend on the choice of $\lambda'$
because of our condition (1') on $\gamma$.

The fact that conditions (1) and (2) hold follows from the construction.
Condition (3) follows from the way in which we build the isomorphism of
\propref{Rtriviality}.

\part{Basic examples}

\section{The universal example: $\GNb$}   \label{univex}

\ssec{} In the category $\Higgs_{\GTb}(\GNb)$ there is a canonical
tautological object. One of the main steps
in the proof of \thmref{second} is to exhibit the
corresponding canonical object in $\Higgs'_{\GTb}(\GNb)$. This is our
goal in this section.

\ssec{} Consider the canonical $T$-bundle $\L_{\B}=G/U$ over $\B=G/B$ and
let us
denote by $\L_{can}$ its pull-back to $\GTb$
under the natural projection $\GTb\to \B$. This will be the first piece in
the data $(\L_{can},\gamma_{can},\beta_{i,can})$.

When we restrict $\L_{can}$ to $G/T\subset \GTb$, it becomes identified
with $G\to G/T$. Hence for every element $\tilde{w}\in N$ that projects to
$w\in W$, we obtain an isomorphism $\gamma_{can}(n):w^*(\L_{can})\simeq
\L_{can}$
over $G/T$, given by right multipliction by $\tilde{w}^{-1}$ on $G$.

However, when exended to the whole of $\GTb$, the above identification is
meromorphic and the configuration of its zeroes and poles is given by a
divisor on $\GTb$ with values in the cocharacter lattice of $T$.

\begin{thm}   \label{descrofpoles}
For a simple reflection $s_i$, the divisor of the above meromorphic
map $s_i^*(\L_{can})\to \L_{can}$ is given by
$-\check\alpha_i(D^{\alpha_i})$.
\end{thm}

The proof will be given in \secref{proofpoles}.

\smallskip

Since $\R^{w_1\cdot w_2}_X\simeq w_2^*(\R^{w_1}_X)\otimes \R^{w_2}_X$,
\thmref{descrofpoles} implies that for any element $w\in W$, the divisor
of zeroes/poles of the above meromorphic map $w^*(\L_{can})\to \L_{can}$
coincides with $\R^w_{\GTb}$. Hence, we obtain the data of
$\gamma_{can}:N\to \Aut_\R(\L_{can})$.

Finally, we have to specify the data of $\beta_{i,can}$ and check
the compatibility conditions. Let us first consider the case when $[G,G]$
is simply connected. As was explained in \secref{simpleversion}, in this
case the data
of $\beta_{i,can}$ can be recovered from $\gamma_{can}$, once
we check that condition (1') holds.

Thus, let $\alpha_i$ be a simple root and let $\lambda$ be a weight
orthogonal to $\check\alpha_i$. It suffices to check condition (1')
at the generic point of $D^{\alpha_i}$. Let $M_i$ be the corresponding
minimal
Levi subgroup. We have a closed embedding $\overline{M_i/T}\subset \GTb$
(cf. \secref{Levicompat})
and its orbit under the $G$-action is the open subset of $\GTb$ equal to
$G/T\cup (D^{\alpha_i}-\underset{\alpha\neq \alpha_i}\cup
(D^\alpha\cap D^{\alpha_i}))$. In particular, it contains a dense subset
of $D^{\alpha_i}$.

Since all our constructions were $G$-equivariant, this implies that
condition (1')
for $\alpha_i$
is equivalent to the corresponding statement for $M_i$. Moreover, we can
replace $M_i$ by
an isogenous group, namely $[M_i,M_i]\times Z(M_i)$. However, in the
latter case
our compatibility condition becomes obvious, as $\lambda$ factors through
$Z(M_i)$.

\medskip

Now, let $G$ be arbitrary. Choose an isogeny $G'\to G$ such that
$[G',G']$ is simply-connected. The varieties $\GTb$ and $\overline{G'/T'}$
are canonically identified and the $T$-bundle $\L_{can}$ is induced
from the $T'$-bundle $\L'_{can}$ under $T'\to T$. Therefore, once we
know the data of $\beta'_{i,can}$ for $\L'_{can}$ that satisfies
the compatibility conditions, it produces the corresponding data
for $\L_{can}$.

\smallskip

Thus, we have constructed a canonical $G$-equivariant object of
$\Higgs'_{\GTb}$ over $\GNb$.

\section{Some simple cases}  \label{examples}

\ssec{The unramified situation} We call a Higgs bundle $(E_G,\sigma)$
unramified if
$\sigma$ maps $E_G$ to $G/N$. Such a map amounts to a reduction of the
structure group from $G$ to $N$.
The category of unramified Higgs bundles is therefore
equivalent to the category of principal $N$-bundles.

The functor $F:\Higgs(X)\to \on{Cam}(X)$ sends an $N$-bundle $E_N$ to
$\tilX:=T\backslash E_N$,
which is a principal $W$-bundle over $X$ (i.e. an \'etale $W$-cover).

In this case the assertion of \thmref{second} is quite evident.

\ssec{$G=SL(2)$.} Fix an $S_2$-cover
$p:\tilX\to X$ and
consider the subsheaf of $p_*(\O_{\tilX})$ consisting of
$S_2$-anti-invariants. We will denote
it by $\z_X$.

It is easy to see that the category
$\Higgs_{\tilX}(X)$ is canonically equivalent to the category
of pairs $(L,\gamma')$, where $L$ is a line bundle on $\tilX$ and
$\gamma'$ is an isomorphism
$\on{det}(p_*(L))\simeq \O_X$.

Let $D_X\subset \tilX$ be the ramification divisor and let
$R_X$ be the corresponding line bundle (cf. \secref{secramification}).
It is easy to see that the category $\Higgs'_{\tilX}(X)$
(which in our case is equivalent to its simplified version
$\Higgs''_{\tilX}(X)$)
consists of pairs $(L,\gamma)$, where $L$ is a line bundle on $\tilX$ and
$\gamma$
is an isomorphism $(-1)^*(L^{-1})\otimes R_X\simeq L$ such that the
composition
$$L\otimes (-1)^*(R_X)\simeq(-1)^*((-1)^*(L^{-1})\otimes R_X)
\overset{(-1)^*(\gamma)}\longrightarrow (-1)^*(L)\overset{\gamma}\simeq
L\otimes (-1)^*(R_X)$$
is minus the identity map.

Let us visualize the equivalence
$\Higgs_{\tilX}(X)\simeq\Higgs'_{\tilX}(X)$ of \thmref{second}
in this case. Indeed, for any line bundle
$L$ on $\tilX$ we have a canonical $S_2$-equivariant isomorphism
$$p^*(\on{det}(p_*(L)))\otimes R_X\simeq L\otimes (-1)^*(L).$$
Therefore, a data of $\gamma'$ defines the data of $\gamma$, and it is
easy to see that this
sets up an equivalence.

\ssec{$G=PGL(2)$.}
In this case the only coroot is non-primitive, so one has to work a little
harder.

By definition, objects of $\Higgs'_{\tilX}(X)$ are the following data:

\begin{itemize}

\item
A line bundle $L$ on $\tilX$.
\item
An $S_2$-equivariant isomorphism of line bundles
$\gamma:L\otimes (-1)^*(L)\simeq R^{\otimes 2}_X$.
\item
An identification $\beta:L|_{D_X}\to R_X|_{D_X}$,
which is compatible in the obvious
sense with the restriction of $\gamma$ to $D_X$.

\end{itemize}

Let us make the statement of \thmref{second} explicit in this case too.
Starting from an object $(E_G,\sigma,t)$
in $\Higgs_{\tilX}(X)$ we can locally choose a
principal $SL(2)$-bundle $E^1_G$, which induces $E_G$. Then
$(E^1_G,\sigma,t)$ is an $SL(2)$-Higgs bundle.
Using the above analysis for $SL(2)$, we can attach to it
a pair $(L^1,\gamma^1)$, where $L^1$ is a
line bundle on $\tilX$ and $\gamma^1:(-1)^*((L^1)^{-1})\otimes R_X\simeq
L^1$.

The corresponding object of $\Higgs'_{\tilX}(X)$ is constructed as
follows.
We define the line bundle $L$ as
$(L^1)^{\otimes 2}$ and $\gamma:=(\gamma^1)^{\otimes 2}$. The data of
$\beta$
comes from the sequence of isomorphisms
$$(L^1)^{-1}\otimes R_X|_{D^X}\simeq
(-1)^*((L^1)^{-1})\otimes R_X|_{D^X}\overset{\gamma^1}\simeq L^1|_{D_X}.$$

If we choose a different lifting of $E_G$ to an $SL(2)$-bundle, the
corresponding
$L^1$ will be modified
by tensoring with $p^*(L^0)$, where $L^0$ is a line bundle on $X$ with
$(L^0)^{\otimes 2}\simeq \O$, which
will not affect the resulting $(L,\gamma,\beta)$.

It is an easy exercise to check that the above construction defines an
equivalence
of categories.

\section{Spectral covers versus cameral covers for $G=GL(n)$}  \label{gln}

\ssec{} Observe first that a regular centralizer in $\gl(n)$ is the same
as an $n$-dimensional
associative and commutative subalgebra in $\on{Mat}(n,n)$ generated by one
element.

\begin{defn}
An $n$-sheeted spectral cover of a scheme $X$ is a finite flat scheme
$p:\Xb\to X$ such that
$p_*(\O_{\Xb})$ has rank $n$ and is locally uni-generated as a sheaf of
algebras.
\end{defn}

Thus, a Higgs bundle for $\gl(n)$ is the same
as a rank-$n$ vector bundle $E$
and an $n$-sheeted spectral cover $\Xb\to X$
with an embedding of bundles of algebras
$p_{*}(\O_{\Xb}) \hookrightarrow \on{End}_{\O_X}(E)$. This is equivalent
to saying that $E$
is a line bundle over $\Xb$.

\medskip

In this section we will analyze the connection of this description of
Higgs bundles
for $GL(n)$ with the one given by \thmref{second}. The starting point is
the observation that
the category of $S_n$-cameral covers of $X$ is naturally equivalent to the
category of $n$-sheeted
spectral covers. Let us describe the functors in both directions:

\smallskip

Given an $S_n$-cameral cover $\tilX\to X$, we define the
scheme $\overline{\tilX}$ as $S_{n-1}\backslash\tilX$.
Conversely, given an $n$-sheeted
spectral cover $\Xb\to X$, we define $\widetilde{\Xb}$ to be
the scheme that represents
the functor of orderings of the sheets of $\Xb\to X$. This functor
attaches to a scheme $S$ the set of data consisting of

\hskip2cm (A map $S\to X$ and $n$ sections $t_i:S\to
\Sbar:=\Xb\underset{X}\times S$), \newline
such that the characteristic polynomial of the multiplication action on
$p_*(\O_{\Xb})$ of any function $f \in \O_{\Sbar}$
equals $\underset{i}\Pi (Y-f\circ t_i)$, where $Y$ is an indeterminate.

\smallskip

It is easy to see that this functor is indeed representable by a scheme
finite over $X$.
The group $S_n$ acts on $\widetilde{\Xb}$ by permuting the $t_i$'s.

\begin{prop}
The functors $\tilX\to \overline{\tilX}$ and
$\Xb\to \widetilde{\Xb}$ send cameral covers to spectral
covers and spectral covers to cameral covers, respectively.
Moreover, they are inverses of one another.
\end{prop}

\begin{proof}

Let us consider first the universal situation:
$X_0=\on{Spec}(\CC[a_0,...,a_{n-1}])$,
$\tilX_0=\on{Spec}(\CC[x_1,...,x_n])$,
where the $x_i$'s satisfy
$$\underset{i}\Pi (Y-x_i)=Y^n+a_{n-1}\cdot Y^{n-1}+...+a_1\cdot Y+a_0,$$
and $\Xb_0=\on{Spec}(\CC[x_1,a_0,...,a_{n-1}])$, where $x_1$ satisfies
$$x_1^n+a_{n-1}\cdot x_1^{n-1}+...+a_1\cdot x_1+a_0=0.$$

The natural maps $\tilX_0\to X_0$ and $\Xb_0\to X_0$ are a cameral and a
spectral
cover, respectively and it is easy to see
that in this case $\overline{\tilX}_0\simeq \Xb_0$ and
$\widetilde{\Xb}_0\simeq \tilX_0$.

This proves the first assertion of the proposition. Indeed, any cameral
(resp., spectral)
cover is locally induced from $\tilX_0$ (resp., $\Xb_0$).

\medskip

For a spectral cover $\Xb$ there is a natural map
$\overline{\widetilde{\Xb}}\to \Xb$, that attaches to a map $S\to
\widetilde{\Xb}$
given by an $n$-tuple $\{t_1,...,t_n\}$ of maps $t_i:S\to \Sbar$ the
composition
$S\overset{t_n}\to \Sbar\to\Xb$. The resulting map
$\overline{\widetilde{\Xb}}\to \Xb$ is an isomorphism,
because this is so in the universal situation, i.e. for
$\overline{\tilX}_0\to X_0$.

Similarly, we have $n$ maps $\tilX\to \overline{\tilX}$ which correspond
to the natural map
$S_n/S_{n-1}\times\tilX\to \Xb$. We claim that
they define an isomorphism $\tilX \to \widetilde{\overline{\tilX}}$.

Indeed, both the fact that these maps satisfy the condition on the
characteristic polynomial
and that the resulting map is an isomorphism follow from the corresponding
facts for
$\tilX_0$.

\end{proof}

\ssec{}

Thus, fixing a spectral cover and fixing an
$S_n$-cameral cover amounts to the same thing. Now, \thmref{second}
implies that the category $\Higgs'_{\tilX}(X)$ is
equivalent to the category of line bundles on the corresponding
spectral cover $\Xb$. We would like to explain how to see this equivalence
explicitly.

We start with the following observation:

\medskip

Let $\tilX\to X$ be an $S_n$-cameral cover and let $\on{Pic}_{\tilX,n}(X)$
be the groupoid
of $S_n$-equivariant line bundles $L$ on $\tilX$ for which the following
condition holds:

\smallskip

\noindent For every reflection $s_{i,j}\in S_n$ the isomorphism
$$s^*_{i,j}(L)\to L$$ is the identity map on the fixed-point set of
$s_{i,j}$ in $\tilX$.

\begin{prop}  \label{picardforgln}
The pull-back functor establishes an equivalence between the category of
line
bundles on $X$ and $\on{Pic}_{\tilX,n}(X)$.
\end{prop}

Let us see first how this proposition implies what we need:

\medskip

The natural map $\tilX\to\Xb$ is itself an $S_{n-1}$-cameral cover. On the
one hand,
by applying the above proposition to this map we obtain that the category
of
line bundles on $\Xb$ is equivalent to $\on{Pic}_{\tilX,n-1}(\Xb)$.

On the other hand, we claim that $\on{Pic}_{\tilX,n-1}(\Xb)$ is equivalent
to
$\Higgs''_{\tilX}(X)$.

Indeed, let us identify the Cartan group of $GL(n)$ with the product of
$n$ copies of $\GG_m$ and let
$\lambda_n:T\to\GG_m$ be the weight corresponding to the last coordinate.
Then a functor
$\Higgs''_{\tilX}(X)\to \on{Pic}_{\tilX,n-1}(\Xb)$ is given by
$(\L,\gamma)\to L:=\lambda_n(\L)$.
It is easy to see that this is indeed an equivalence.

\ssec{} Now let us prove \propref{picardforgln}. The argument will be a
prototype of the one
we are going to use to prove \thmref{second}.

Given an object $L\in \on{Pic}_{\tilX,n}(X)$ and a point $x\in X$ we must
find an \'etale
neighbourhood of $x$ such that, when restricted to the preimage of this
neighbourhood, $L$ becomes
isomorphic to the unit object in $\on{Pic}_{\tilX,n}(X)$ (i.e. the one for
which
$L=\O_{\tilX}$ with the tautological $S_n$-structure).

First, it is easy to reduce the statement to the case when the
ramification
over $x$ is the maximal possible,
i.e. when $x$ has only one geometric preimage $\tilde{x}$ in $\tilX$.
Further, we can assume that $X$ (and therefore
also $\tilX)$ is a spectrum of a local ring.

Choose some trivialization of $L$. Its discrepancy with the
$S_n$-equivariant structure
is a $1$-cocycle
$S_n\to \Hom(\tilX,\GG_m)$. We must show that this cocycle is homologous
to $0$.

\smallskip

Let $K$ denote the kernel of the map $\Hom(\tilX,\GG_m)\to\GG_m$
given by the evaluation at $\tilde{x}$.
Our condition on $L$ implies that the above cocycle
$S_n\to \Hom(\tilX,\GG_m)$ takes
values in $K$. However, since $\tilX$
is local, $K$ is divisible and torsion-free. Hence $H^1(S_n,K)=0$, so our
cocycle is cohomologically trivial.

\part{Basic structure results over $\GNb$}

\section{The structure of $\GNb$}   \label{structure}

\ssec{} The rest of the paper is devoted to the proofs of various results
announced in the
previous sections. We start with the proof of \propref{GmodN}.

\begin{proof}

First we need to show that the map $\phi:\gf_{reg}\to \on{Ab}^r$ is
well-defined, which is equivalent
to saying that the projection $\Gamma_{reg}\to\gf_{reg}$ is an
isomorphism.

Since the latter projection is proper and $\gf_{reg}$ is reduced, it is
enough to show
that the scheme-theoretic preimage in $\Gamma_{reg}$ of every
$x\in\gf_{reg}$
is isomorphic to $\Spec(\CC)$.

\medskip

This is clear on the level of $\CC$ points, since by definition of regular
elements, the only abelian $r$-dimensional subalgebra in $\gf$ that
contains $x$ is its own centralizer.

\medskip

For $\af\in \on{Ab}^r$, the tangent space $T_\af(\on{Ab}^r)$ can be
identified with the space of maps $T:\af\to\gf/\af$ that satisfy:
$$\forall y_1,y_2\in\af,\,\,[T(y_1),y_2]+[y_1,T(y_2)]=0\in\gf.$$

\smallskip

We claim that the tangent space to $\Gamma_{reg}\cap (\on{Ab}^r\times
x)$ at $\af\times x$ is zero. Indeed, this is the space
of maps $T:\af\to\gf/\af$ as above, for which, moreover
$[T(y),x]=0,\,\,\forall y\in\af$. However,
since $\af=Z_\gf(x)$, any such $T$ is identically zero.

\smallskip

This implies that $T_{\af\times x}(\Gamma_{reg}\cap \on{Ab}^r\times
x)=0$ which means that $\Gamma_{reg}\cap (\on{Ab}^r\times x)$
is reduced, i.e. $\simeq\Spec(\CC)$.

\medskip

Now let us show that $\phi$ is smooth. Let $\af\in \on{Ab}^r$ be equal to
$\phi(x)$.
Using the above description of the tangent space to $\on{Ab}^r$
it is easy to see that $d\phi$ sends an element $u\in\gf\simeq
T_x(\gf_{reg})$ to the unique map
$T:\af\to\gf/\af$ that satisfies:  $$[T(y),x]+[y,u]=0,\,\,\forall
y\in\af.$$

\smallskip

Consider now the map
$\ev: T_\af(\on{Ab}^r)\to \gf/\af$ given by
$T\to T(x)$. The above description of $d\phi$ implies that the composition
$$\gf\simeq T_x(\gf_{reg}){\overset{d\phi}\longrightarrow}
T_\af(\on{Ab}^r){\overset{\ev}\longrightarrow} \gf/\af$$
coincides with the tautological projection $\gf\to\gf/\af$.

\smallskip

However, since $x$ is regular, the fact that
$[T(x),y]=-[x,T(y)],\,\,\forall y\in\af$
implies that $\ev$ is an injection. We conclude that
$\ev$ is an isomorphism, hence $\on{Im}(\phi)$ is contained in the
smooth locus of $\on{Ab}^r$. Furthermore,
$d\phi$ is surjective, so $\phi$ is smooth as claimed.

\end{proof}

\ssec{} Let $\tilgf$ be the closed sub-variety in $\gf\times\B$ defined by
the
condition:  $(x,\bof')\in\tilgf \text{ if } x\in \bof'$.  Let
$\tilgf_{reg}$
denote the intersection $\tilgf\cap (\gf_{reg}\times\B)$ and let $\tilpi$
denote the projection $\tilgf_{reg}\to\gf_{reg}$. It is clear that as a
variety,
$\tilgf_{reg}$ is smooth and connected, since it is an open
subset in a vector bundle over $\B$.

\begin{prop} \label{GmodT}
There exists a natural $G$-invariant map
$\tilphi:\tilgf_{\on{reg}}\to\GTb$,
such that the following square is Cartesian:
$$
\CD
\tilgf_{reg} @>{\tilphi}>> \GTb \\
 @V{\tilpi}VV @V{\pi}VV  \\
\gf_{reg}  @>{\phi}>> \GNb
\endCD
$$
\end{prop}

\begin{proof}

Consider the fibered product $\GTb{\underset{\GNb}\times}\gf_{reg}$.
By definition of $\GTb$, there is a closed embedding
$$\GTb{\underset{\GNb}\times}\gf_{reg}\to\tilgf_{reg}$$ that sends a
triple $(\af\in\GNb,\bof'\in\B,x\in\gf_{reg})\in
\GTb{\underset{\GNb}\times}\gf_{reg}$ to $(x,\bof')\in\tilgf_{reg}$.

We claim that this embedding is in fact an isomorphism.
Indeed, the statement is obvious over the preimage in $\tilgf_{reg}$ of
the regular semisimple locus of $\gf$.  Therefore, the two schemes
coincide at the generic point of $\tilgf_{reg}$. This implies what we
need,
since $\tilgf_{reg}$ is reduced.

\end{proof}

Now we are ready to prove \propref{WactsonGmodT}:

\begin{proof}

The map $\tilphi:\tilgf_{reg}\to\GTb$ is smooth, since it is a base change
of a smooth map.
Hence, the fact that $\tilgf_{reg}$ is smooth and connected implies that
$\GTb$
has the same properties.

A well-known theorem of Kostant (cf.  \cite{Ko} or \cite{Di}, p. 277) says
that the restriction of the Chevalley map $\gf\to \tof/W$ to $\gf_{reg}$
is
smooth and that it gives rise to a Cartesian square:

$$
\CD
\tilgf_{reg} @>>> \tof  \\
@VVV        @VVV      \\
\gf_{reg} @>>>  \tof/W
\endCD
$$

Therefore, the natural action of $W$ on the preimage in $\tilgf$ of the
regular semisimple locus in $\gf$ extends to the whole of $\tilgf_{reg}$.
The same is true for $\GTb$, because the map $\tilphi$ is flat and
surjective.
The \'etale local isomorphism follows from comparison of our Cartesian
square with
that of \propref{GmodT}.

\end{proof}

\ssec{}

Now let us prove \propref{orbits}.

\begin{proof}

Let $\Delta'$ be as in the formulation of the proposition. Consider
an element $t\in\tof$ such
that $\alpha(t)=0$ for $\alpha\in\Delta'$ and $\beta(t)\neq 0$
for $\beta\notin \Delta'$.

In this case $\mf:=Z_{\gf}(t)$ is a Levi subalgebra of $\gf$.
Let $M$ be the corresponding Levi subgroup. It is well-known
that $\mf\cap\bof$ is a Borel subalgebra in $\mf$. Let $u$
be an element in the unipotent radical of $\mf\cap\bof$, which
is regular with respect to $M$.

We then see that $x=t+u$ is a regular element in $\gf$, since
$Z_{\gf}(x)=Z_{\mf}(u)$. It is known that if a Borel subalgebra
contains a regular element, then it also contains its centralizer (cf.
\lemref{inclusion}).
Therefore, $(Z_{\mf}(u),\bof)\in\GTb$. Moreover,
it is easy to see that every pair $(\af,\bof')\in\GTb$ is $G$-conjugate to
one of the above form.

To conclude the proof, it remains to show that
$(Z_{\mf}(u),\bof)\in \underset{\alpha\in\Delta'}\cap (D^\alpha)\setminus
\underset{\beta\notin\Delta'}\cup (D^\beta)$. For that, it suffices to
show that
the image of $(t+u,\bof)$ as above under $\tilgf_{reg} \to \tof$
belongs to the corresponding locus of $\tof$. However, the above
image is just $t$, which makes the assertion obvious.

\end{proof}

\ssec{Levi subgroups}   \label{Levicompat}
Let $J\subset I$ be a subset. It defines a root subsystem $\Delta_J$
and let $M_J$ (resp., $P_J\subset G$, $W_J\subset W$) denote the
corresponding standard Levi
subgroup (resp., standard parabolic, Weyl subgroup). Let $N_{M_J}$ be
the intersection $M_J\cap N$, which is the normalizer of $T$ in $M_J$.

\smallskip

It is easy to see that the natural map $M_J/N_{M_J}\to G/N$ extends to a
map
$i_J:\overline{M_J/N_{M_J}}\to \GNb$. In fact, $\overline{M_J/N_{M_J}}$ is
a
closed sub-variety of $\GNb$ which corresponds to
$\{\af\in\GNb\,|\,\af\subset\mf_J\}$.

\begin{prop}  \label{pullbackGmodT}
There is a canonical $W$-equivariant isomorphism:
$$\tili_J:W{\overset{W_J}\times}\overline{M_J/T}\simeq\overline{M_J/N_{M_J}}
{\underset{\GNb}\times}\GTb.$$
\end{prop}

\begin{proof}

First, we have a natural closed embedding
$$\overline{M_J/T}\to
\overline{M_J/N_{M_J}}{\underset{\GNb}\times}\GTb\subset \GTb.$$
Its image consists of pairs $(\af,\bof')\in\GTb$ such that
$\af\subset \mf_J$ and $\bof'\subset \pf_J:=\on{Lie}(P_J)$.

This map is compatible with the $W_J$-action. Hence, it extends
to a finite map
$$\tili_J:W{\overset{W_J}\times}\overline{M_J/T}\to\overline{M_J/N_{M_J}}
{\underset{\GNb}\times}\GTb.$$

Since both varieties are smooth, in order to prove that $\tili_J$ is an
isomorphism, it suffices to do so over the open part, i.e. over
$M_J/N_{M_J}$. However, in the latter case, the assertion becomes obvious.

\end{proof}

It is easy to see that the $G$-orbit of
$\overline{M_J/N_{M_J}}{\underset{\GNb}\times}\GTb\subset\GTb$
(resp., $\overline{M_J/T}\subset \GTb$) is
the union of those $D^{\Delta'}$ for which $\Delta'$
is $W$-conjugate to a subset of $\Delta_J$ (resp.,
$\Delta'\subset\Delta_J$.)

\section{The group-scheme of centralizers}  \label{grschcentr}

In this section we will formulate two
basic theorems, \thmref{centralizers} and \thmref{conjugation},
which will be used for the proof of our first main result, \thmref{first}.

\ssec{The universal centralizers $\Z$ and $\z$} \label{centralizer_families}
Consider the constant group-scheme $G\times\GNb$ over $\GNb$, and let
$\Z\subset G\times\GNb$ be its closed group-subscheme of ``centralizers''.
In other words, $\Z$ is defined by the condition that
$(g\in G,\af\in\GNb)\in \Z$ if $g$ commutes with $\af$.
Clearly, $\Z$ is equivariant with respect to the $G$-action on $\GNb$.

\smallskip

Note that the corresponding bundle $\z$ of Lie algebras can be identified
with the tautological rank $r$ vector bundle over $\GNb$ which comes
from the embedding $\GNb\subset Gr^r_\gf$. Another interpretation of this
$\z$, considered as a subbundle of the trivial bundle
$\gf \times \GNb$, is that it is the family $\z_{\GNb}$ of
centralizers of the universal Higgs bundle on $\GNb$,
which was studied in detail in
\secref{univex}. (Recall from \secref{intrHiggs} that a Higgs bundle
$(E_G,\sigma)$ on any $X$ determines, and is determined by, a subbundle
$\z_{X}$ consisting of regular centralizer subalgebras of the adjoint bundle
$\gf_{E_G}$.)

\begin{prop} \label{centrproperties}
The group-scheme $\Z$ is commutative and smooth over $\GNb$ and is
irreducible as a variety.
\end{prop}

\begin{proof}

Let $\Z'$ be the group-subscheme of $G\times\gf_{reg}$ over
$\gf_{reg}$ defined by the condition:
$$\Z':=\{(g,x)\in G \times \gf_{reg} \,|\,\on{Ad}_g(x)=x\}.$$
First, let us show that $\Z'$ is commutative and smooth over $\gf_{reg}$.

Let $(g,x)$ be a $\CC$-point of $\Z'$. The tangent space to
$\Z'$ at $(g,x)$ consists of pairs
$(\xi,y)\in\gf\times\gf$ such that $\on{Ad}_g([x,\xi])=\on{Ad}_g(y)-y$.
The differential of the map $\Z'\to\gf_{reg}$ sends
$(\xi,y)$ to $y$. We claim that it is surjective.

It is known that if $G$ is of adjoint type, then the centralizer of
every regular element is connected. (In particular, each $Z_G(x)$ is
commutative; this holds even if $G$ is not of adjoint type.). Therefore,
$\underset{g\in
Z_G(x)}{\on{Span}}(\on{Ad}_g(y)-y)=\on{Im}(\on{ad}_{Z_{\gf}(x)})$.
However, the latter, as we saw in the proof of \propref{GmodN}, coincides
with $\on{Im}(\on{ad}_x)$,
since $x$ is regular.

To prove that $\Z'$ is smooth over $\gf_{reg}$, it remains to observe that
the fibers of $\Z'$ are smooth (since they are algebraic groups in
char.$0$)
and all have dimension $r$, by the definition of $\gf_{reg}$. The fact that
$\Z'$
is commutative was established in the course of the above argument.

\medskip

Now let us prove the assertion for $\Z$. We have a natural closed
embedding
$\Z\underset{\GNb}\times\gf_{reg}\to\Z'$,
which is an isomorphism over the regular semisimple locus of $\gf_{reg}$.
Hence,
it is an isomorphism, because $\Z'$ is reduced. Therefore, since the map
$\phi:\gf_{reg}\to\GNb$
is flat and surjective, this shows that $\Z$ is commutative and
smooth over $\GNb$. It is
irreducible, because this is obviously true over $G/N$.

\end{proof}

\ssec{The group scheme $\TT$.} \label{TT}
Now we will introduce another group-scheme
over $\GNb$, seemingly of a different
nature. Recall the sheaves $\Tb_{\tilX},T_{\tilX}$ introduced in
section \ref{HiggsIsGerbe}.

\smallskip

Consider the contravariant functor Schemes $\Rightarrow$ Groups which
assigns to a scheme $S$ the set of pairs

\smallskip

\hskip2cm (A map $S\to\GNb$, a $W$-equivariant map
$\tilS:=S{\underset{\GNb}\times}\GTb\to T).$

\smallskip

It is easy to see that this functor is representable by an abelian
group-scheme over $\GNb$, which we will denote by $\overline{\TT}$.
Therefore, once $S\to \GNb$ is fixed,
$\Hom_{\GNb}(S,\overline{\TT})\simeq \Gamma(S,\overline{T}_{\tilS})$.
In other words, $\overline{\TT}$ represents the sheaf $\overline{T}_{\GTb}$
on $\on{Sch}_{et}(\GNb)$.

Clearly, the $G$-action on $\GTb$ gives rise to a $G$-action on
$\overline{\TT}$.

\medskip

We define the open group-subscheme $\TT$ of $\overline{\TT}$ by the
following condition (**):

\smallskip

\noindent $\Hom(S,\TT)$ consists of those pairs
$(S\to\GNb,\tilS\to T)$ as above,
for which for every root $\alpha$ the composition
$$S\underset{\GNb}\times D^\alpha\hookrightarrow
\tilS\to T\overset{\alpha}\to {\mathbb G}_m$$
avoids $-1\in {\mathbb G}_m$.

\smallskip

Since for any map $S\to \TT$, the above composition takes values
in $\pm 1\subset \GG_m$,
condition (**) is equivalent to condition (*)
in the definition of the sheaf $T_{\tilS}$ (cf. \secref{sheaf}):
for a fixed map $S\to\GNb$, $\Hom_{\GNb}(S,\TT)\simeq \Gamma(S,T_{\tilS})$,
i.e., the group-scheme $\TT$ represents the sheaf $T_{\GTb}$ on $\on{Sch}_{et}(\GNb)$.

\ssec{} A remarkable fact is that the group-schemes $\Z$ and $\TT$ are
canonically isomorphic.
Here we will construct a map between them in one direction.

\smallskip

Let $\Bu$ denote the universal group-scheme of
Borel subgroups over $\B$. Let us denote by
$\widetilde{\Bu}$ its pull-back to $\GTb$.
In addition, let us denote by $\tilZ$ the pull-back of
$\Z$ to $\GTb$.

Both $\widetilde{\Bu}$ and $\tilZ$ are group-subschemes
of the constant group-scheme
$G\times\GTb$.

\begin{lem} \label{inclusion}
$\tilZ$ is a closed group-subscheme of $\widetilde{\Bu}$.
\end{lem}

Indeed, since $\tilZ$ is reduced and irreducible, it suffices to check
that
over $G/N$, $\tilZ$ is contained in $\widetilde{\Bu}$. However, this
is obvious.

We have a natural projection $\Bu\to T\times\B$. By composing it with the
inclusion of \lemref{inclusion}, we obtain a map
$$\Z\underset{\GNb}\times\GTb\to T.$$

This map respects the group law on $\Z$ and $T$
and commutes with the $W$-action.
(This is because it suffices to check both facts
after the restriction to $G/N$, where they become obvious.)

\smallskip

Hence, we obtain a homomorphism of group-schemes $\overline{\chi}:\Z\to \overline{\TT}$.

\begin{thm} \label{centralizers}
The above map $\overline{\chi}:\Z\to\overline{\TT}$ defines an isomorphism
${\chi}:\Z\to \TT$.
\end{thm}

The proof will be given in the next section.

\ssec{} Now we will formulate the second key result which will be used in
the proof of \thmref{first}.

\smallskip

Consider the functor that assigns to a scheme $S$ the set of triples
$(\GNb^1_S,\GNb^2_S,\nu)$, where $\GNb^1_S$ and $\GNb^2_S$ are two
$S$-points of $\GNb$ and $\nu$ is a $W$-equivariant isomorphism
$$\nu:\tilS^1\to\tilS^2,$$
where $\tilS^i$ is the $W$-cover of $S$ induced
by $\GNb^i_S$ from $\pi:\GTb\to\GNb$.

\smallskip

It is easy to see that this functor is representable.
Let $\H$ denote the representing scheme.
Since the $W$-cover $\GTb\to\GNb$ is $G$-equivariant, we obtain a natural
map
$\xi:G\times \GNb\to \H$ which covers the map $G\times \GNb
\overset{\on{Action}\times\on{id}}\longrightarrow\GNb\times\GNb$.

\begin{thm} \label{conjugation}
The above map $\xi:G\times\GNb\to \H$ is smooth and surjective.
\end{thm}

This theorem will be proven in \secref{proofconj}.

\ssec{} \label{secstabilizers}

The scheme $\H$ lives over $\GNb\times\GNb$. Let $\H_\Delta$
denote its restriction to the diagonal. By definition,
$\H_\Delta$ is a group-scheme
over $\GNb$ which represents the functor of $W$-equivariant automorphisms
of
$\GTb$ over $\GNb$.

Let $\on{St}\subset G\times\GNb$ be the closed group-subscheme of
stabilizers, i.e.
$$(g,\af)\in \on{St} \text{ if } \on{Ad}_g(\af)=\af.$$
Obviously, $\Z$ is a closed normal group-subscheme of $\on{St}$.

The map $\xi:G\times\GNb\to\H$ gives rise to a map $\xi_\Delta:\on{St}\to
\H_\Delta$.

\begin{prop} \label{stabilizer}
$\H_\Delta$ represents the quotient group-scheme $\on{St}/\Z$.
\end{prop}

\begin{proof}

\thmref{conjugation} implies that the map $\xi_\Delta:\on{St}\to
\H_\Delta$
is smooth and surjective. Therefore, all we have to show is that
if $S\to \on{St}$ is a map such that the induced automorphism of $\tilS$
is trivial, then $S$ maps to $\Z$.

Observe that $\H_\Delta$ acts on $\TT$ via its action on $\GTb$.
Since the isomorphism $\chi:\Z\to \TT$ is $G$-equivariant, we obtain
a commutative diagram of actions:
$$
\CD
\on{St}\protect\underset{\GNb}\times \Z  @>>> \Z \\
@V{\xi_\Delta\times\chi}VV          @V{\chi}VV  \\
\H_\Delta\underset{\GNb}\times \TT  @>>> \TT,
\endCD
$$
where the top horizontal arrow is the adjoint action.

\smallskip

Therefore, if a map $S\to \on{St}$ induces the trivial automorphism of
$\tilS$,
its adjoint action on $\Z$ is trivial too. But this means that it factors
through
$\Z$.

\end{proof}

Similarly, one shows:

\begin{cor} \label{quotient}
The scheme $\H$ represents the quotient
group-scheme $(G\times \GNb)/\Z$.
\end{cor}

\ssec{}

Here is one more interpretation of \thmref{conjugation}:

Clearly, the scheme $\H$ with its two projections to $\GNb$
is a groupoid over that latter scheme. According to
\thmref{conjugation}, the above projections are smooth
and, therefore, we can consider the algebraic stack $\H\backslash (\GNb)$.

\begin{cor}
The stack $\H\backslash (\GNb)$ is canonically isomorphic to
the stack $\bfCam$ of \secref{intr covers}.
\end{cor}

\section{Proof of \thmref{centralizers}}

\ssec{}
We start by establishing a result on compatibility of our objects with
restrictions to Levi subgroups. We then verify the Lie-algebraic version
of the theorem by restricting to an $\sl(2)$ subalgebra, and finally we
refine this to prove the desired group-theoretic version.

\ssec{}

Let $M=M_J$ be a standard Levi subgroup of $G$ (cf. \secref{Levicompat})
and
let $\Z_M$ be the corresponding sheaf of centralizers
over $\overline{M/N_M}$.

On the one hand, there is a natural closed embedding
$$\Z_M\hookrightarrow
i^*_J(\Z):=\overline{M/N_M}{\underset{\GNb}\times}\Z.$$
On the other hand, we have the group scheme $\overline{\TT}$ over
$\GNb$, as well as the group-scheme
$\overline{\TT}_M$ over $\overline{M/N_M}$. This time,
by \propref{pullbackGmodT}, we have a canonical isomorphism
$$\overline{\TT}_M\simeq
i^*_J(\overline{\TT}):=\overline{M/N_M}{\underset{\GNb}\times}\overline{\TT}.$$
Moreover, it induces an isomorphism $\TT_M\simeq i^*_J(\TT)$, since if a
root
$\alpha$ is not $W$-conjugate to a root in $M$, then $s_\alpha$
has no fixed points on $W{\overset{W_J}\times}\overline{M/T}$.

\begin{prop}  \label{reductiontoLevi}

The map $\Z_M \to i^*_J(\Z)$
is an isomorphism. Moreover, the diagram
$$
\CD
\Z_M @>{\chi_M}>>  \overline{\TT}_M  \\
@VVV          @VVV          \\
i^*_J(\Z)  @>{i^*_J(\chi)}>> i^*_J(\overline{\TT})
\endCD
$$
is commutative.

\end{prop}

\begin{proof}

The map $\Z_M\hookrightarrow i^*_J(\Z)$ is an isomorphism
because it is a closed embedding and at the same time
an isomorphism over the generic point of $\overline{M/N_M}$. Commutativity
of the diagram can be checked over the preimage of $M/N_M$, in which case
it becomes obvious.

\end{proof}

\ssec{} We will now prove the assertion of
\thmref{centralizers} on the Lie-algebra level.

\smallskip

Let $\tt$ denote the sheaf of Lie algebras corresponding
to $\overline{\TT}$. Obviously, it is isomorphic to $\on{Lie}(\TT)$
as well. By definition, we have: $\tt\simeq
(\tof\otimes\pi_*(\O_{\GTb}))^W$.
Since $\pi_*(\O(\GTb))$ is locally isomorphic to
$\O_{\GNb}\otimes \CC[W]$, $\tt$
is a vector bundle of rank $r$ over $\GNb$.

On the other hand, recall that in subsection \ref{centralizer_families}
we defined the sheaf $\z$ of Lie algebras corresponding to $\Z$. Our map
$\overline{\chi}:\Z\to\overline{\TT}$ induces a map
$d{\overline{\chi}}:\z\to\tt$ which, for simplicity, we abbreviate as
$d\chi:\z\to\tt$.

\begin{prop} \label{Liealgebrasfirst}
The map $$d\chi:\z\to\tt$$ is an isomorphism.
\end{prop}

\begin{proof} The proof will consist of two steps. The first step will be
a reduction to the case of $SL(2)$ and the second one will be a proof of
the assertion for $SL(2)$.

\smallskip

\noindent {\bf Step 1.}

Both $\z$ and $\tt$ are vector bundles of rank $r$ over $\GNb$ and
the map $d\chi$ is clearly an isomorphism over $G/N$. Since the variety
$\GNb$ is smooth, it remains to show that $d\chi$ is an isomorphism on an
open subset of $\GNb$ whose complement has codimension at least $2$.

\smallskip

It follows from \secref{Levicompat} that such an open subset
is formed by the union of the $G$-orbits of the images of
$i_J(\overline{M_J/N_{M_J}})$, where $J=\{\alpha_j\}$ for all
simple roots $\alpha_j$.

Therefore, by $G$-equivariance and by \propref{reductiontoLevi},
it suffices to show that
the map $$d\chi_{M_J}:\z_{M_J}\to\tt_{M_J}$$ is an isomorphism.
This reduces us to the case when $G$ is a reductive group of semi-simple
rank $1$.

\smallskip

Moreover, the statement is clearly invariant under isogenies, so we may
replace $G$ by $Z^0(G)\times [G,G]$. Clearly, the assertion in such a
case is equivalent to the one for $[G,G]$, which in
turn can be replaced by $SL(2)$.

\smallskip

\noindent {\bf Step 2.}

For $G=SL(2)$, the variety $\GNb$ can be identified with $\PP^2$
in such a way that the
sheaf $\z$ goes over to $\O(-1)$. Moreover,
$\GTb\to\GNb$ can be identified with the $S_2$-cover
$\pi:\PP^1\times\PP^1\to\PP^2$.

\smallskip

To prove the assertion, it is enough to show that $\tt$ has degree $-1$,
since any non-zero map between two line bundles of the same degree is
automatically an isomorphism.

\smallskip

By definition, $\tt$ is the $\O(\PP^2)$-module of anti-invariants of
$S_2$ in $\pi_*(\O(\PP^1\times\PP^1))$.  Therefore,
$$\tt\simeq\det(\pi_*(\O(\PP^1\times\PP^1)))\simeq \O(-1).$$

\end{proof}

\ssec{}
Now we will check that the map $\chi$ induces an isomorphism between
$\CC$-points of $\Z$ and $\TT$.
Evidently, this assertion, combined with
\propref{Liealgebrasfirst} and \propref{centrproperties}, implies
\thmref{centralizers}.

\medskip

Let $\af\in\GNb$ be the centralizer of a regular element $x\in\gf$.
As we saw in the proof of \propref{centrproperties}, on the one hand,
the fiber of $\Z$ at $\af=Z_\gf(x)$ can be
identified with $Z_G(x)$. On the other hand, the fiber of $\TT$ at
$\af$ can be identified with $\Hom_W(\B^x,T)^{**}$, where $\B^x$ is the
fixed
point scheme of the vector field induced by $x$ on $\B$ and the
super-script $**$ corresponds to the (**) condition in the definition
of $\TT$.

\smallskip

Let $x=x^{ss}+x^{nil}$ be the Jordan decomposition of $x$.
We can assume that $Z_{\gf}(x^{ss})$
is a
standard Levi subalgebra $\mf$ and $x^{nil}$ is a regular nilpotent
element
in $\mf$.
Using \propref{reductiontoLevi}, we can replace $G$ by $M$ and hence we
can
assume that $x^{ss}$ is a central element in $\gf$.

\smallskip

There are natural embeddings $Z(G)\times\GNb\to \Z$ and
$Z(G)\times\GNb\to\TT$, which make the diagram
$$
\CD
Z(G)\times\GNb  @>>> \Z \\
@V{\on{id}}VV   @V{\chi}VV \\
Z(G)\times\GNb @>>> \overline{\TT}
\endCD
$$
commute.

\smallskip

\begin{prop}  \label{jordan}
Let $x$ be a regular nilpotent element and let
\begin{align*}
&Z_G(x)=Z_G(x)^{ss}\times Z_G(x)^{nil}  \\
&\Hom_W(\B^x,T)^{**}=
\Hom_W(\B^x,T)^{ss,**}\times \Hom_W(\B^x,T)^{nil,**}
\end{align*}
be the Jordan decompositions of the fibers of $\Z$ and $\TT$ at
$Z_{\gf}(x)$. Then the embedding of $Z(G)$ induces isomorphisms:
$$Z(G)\simeq Z_G(x)^{ss}\text{ and } Z(G)\simeq\Hom_W(\B^x,T)^{ss,**}.$$
\end{prop}

It is clear, first of all, that this proposition implies the theorem.
Indeed, it is enough to show that $\chi$ induces an isomorphism
$Z_G(x)^{nil}\to \on{Hom}_W(\B^x,T)^{nil,**}$. But since these groups are
unipotent, our assertion follows from the corresponding assertion on the
Lie-algebra level, which has been proven before.

\begin{proof}
The fact that $Z(G)\simeq Z_G(x)^{ss}$ is an immediate consequence of the
fact that in a group of adjoint type centralizers of regular elements are
connected.

\smallskip

To prove that $Z(G)\simeq\Hom_W(\B^x,T)^{ss,**}$, let us observe that
if $x$ is a regular nilpotent element, $\B^x$ is a local non-reduced
scheme.
Its closed point, viewed as a point of $\GTb$, belongs to the intersection
of all the $D^\alpha$'s.

Let $\Hom_W(\B^x,T)_1$ be the sub-group of $\Hom_W(\B^x,T)$ which
corresponds
to maps $\B^x\to T$ that send the closed point of $\B^x$ to the identity
in $T$. Clearly, $\Hom_W(\B^x,T)_1$ is unipotent and
$\Hom_W(\B^x,T)\simeq \Hom_W(\B^x,T)_1\times T^W$ is the Jordan
decomposition of $\Hom_W(\B^x,T)$.

\smallskip

The proof is concluded by the observation that
$$Z(G)=\{t\in T^W\, \ | \ \, \alpha(t)=1,\,\forall\alpha\in\Delta\},$$
which is exactly the (**) condition.

\end{proof}

\section{Proof of \thmref{conjugation}} \label{proofconj}

\ssec{} We will need an additional property of the isomorphism $\chi$.

\medskip

By definition, we have a canonical $W$-equivariant map
$$\tt\underset{\GNb}\times\GTb \to \tf,$$
hence we obtain a map $\tt\to \tf/W$.

\begin{lem} \label{Chevy}
The above map coincides with the composition
$$\tt\overset{\chi^{-1}}\longrightarrow \z\hookrightarrow \gf\times\GNb\to
\gf\to
\tf/W,$$
where the last arrow is the Chevalley map.
\end{lem}

The proof follows from the fact that the two maps coincide
over $G/N$.

\ssec{}

Since $G\times \GNb$ is smooth,
to prove the theorem, we need to show that any map
$S\to \H$ can be lifted, locally in the \'etale topology,
to a map $S\to G\times \GNb$.

Thus, let $\af^1$ and $\af^2$ be two $S$-points of $\GNb$
and let $\nu:\tilS^1\to\tilS^2$ be an isomorphism
between the corresponding cameral covers. The maps $\af^i$ give rise to
vector subbundles
$\z^i_S\subset \gf\otimes \O_S$, and \thmref{centralizers} implies that
$$\z^i_S\simeq \Hom_{W,\O_S}(\tf^*,\O_{\tilS^i}), \,\, i=1,2.$$
Therefore, the data of $\nu$ defines an isomorphism of vector bundles
$\nu':\z^1_S\to \z^2_S$.

\medskip

By \propref{GmodN} we can find a section $x_S^1\in \z^1_S$, such that
$\z^1_S=Z_{\gf}(x^1_S)$.
Let $x^2_S\in \z^2_S$ be the image of $x_S^1$ under $\nu'$. By making the
choice of $x_S^1$ sufficiently
generic, we can assume that $x^2_S$ is regular, i.e. that
$\z^2_S=Z_{\gf}(x^2_S)$.

\smallskip

Consider $x^i_S,\, i=1,2$ as maps $S\to \gf_{reg}$. \lemref{Chevy} implies
that their compositions
with the Chevalley map
$$S\overset{x^i_S}\longrightarrow \gf_{reg}\to \tf/W$$
coincide. Now, we have the following general assertion that follows from
smoothness of the Chevalley map
restricted to $\gf_{reg}$:

\begin{lem}
The adjoint action map $G\times \gf_{reg}\to
\gf_{reg}\underset{\tf/W}\times \gf_{reg}$
is smooth and surjective.
\end{lem}

Therefore, locally there exists a map $g_S:S\to G$ that conjugates $x_S^1$
to $x^2_S$.
Then this map conjugates $\z^1_S$ to $\z^2_S$, which is what we had to
prove.

\ssec{Complements}

We conclude this section by two remarks regarding the assertions of
\thmref{centralizers} and
\thmref{conjugation}.

\medskip

First, let us fix a $\CC$-point $\af\in\GNb$ and let $\varphi:\B^{\af}\to
\tf$ be a
$W$-equivariant map, which according to \thmref{centralizers},
is the same as an element $x_{\varphi}\in\af=\z_{\af}$. One may
wonder: how can one express the condition that $x_{\varphi}$
is a regular element of $\af$ in terms of $\varphi$?

\begin{lem}  \label{regcriterion}
The necessary and sufficient condition for $x_{\varphi}$
to be a regular element of $\af$ is that
$\varphi:\B^{\af}\to \tf$ is a scheme-theoretic embedding.
\end{lem}

\begin{proof}

First, one easily reduces the assertion to the case when $\af$ is the
centralizer of a regular nilpotent element, which we will assume.

In this case, $\af$ entirely consists of nilpotents elements.
Let $\on{St}_G(\af)$ be the normalizer of $\af$. Since 
the nilpotent locus in $\gf_{reg}$
is a single $G$-orbit, we obtain that $\af\cap \gf_{reg}$ is
a single $\on{St}_G(\af)$-orbit.

Thus, let $\phi$ be an embedding. To show that $x_{\varphi}$ is regular,
it is enough to show that its centralizer in $\on{St}_G(\af)$ coincides with
$\Z_{\af}$.

By \propref{stabilizer}, the quotient $\on{St}_G(\af)/\Z_{\af}$ maps
isomorphically to the group of $W$-equivariant automorphisms of $\B^{\af}$.
If for some $n\in \on{St}_G(\af)$ we have $\on{Ad}_n(x_{\varphi})=x_{\varphi}$,
then $n$ acts trivially on $\B^{\af}$, since $\phi$ is an embedding.
Hence, $n\in \Z_{\af}$.

\medskip

To prove the implication in the other direction, let us observe that
$\af\cap \gf_{reg}$ is the only $\on{St}_G(\af)$-invariant open
subset of $\af$ consisting of regular elements only.
However, the locus of $\varphi$ that are embeddings is clearly such a
subset.

\end{proof}

Secondly, let us see how \corref{quotient} is related to \propref{orbits}:

\smallskip

Let $\af_1$ and $\af_2$ be two $\CC$-points of $\GNb$. \corref{quotient}
says that
they are $G$-conjugate if and only if $\pi^{-1}(\af_1)\simeq
\pi^{-1}(\af_2)$ as $W$-schemes.
The condition of \propref{orbits} is seemingly weaker (but in fact,
equivalent):
it implies that $\af_1$ and $\af_2$
are $G$-conjugate if and only if $(\pi^{-1}(\af_1))_{red}\simeq
(\pi^{-1}(\af_2))_{red}$
as $W$-schemes.

\part{Proofs of the main results}

\section{Proof of \thmref{first}}

\ssec{} We are going to deduce our theorem from \thmref{centralizers} and
\thmref{conjugation}
combined with the following ``abstract nonsense'' observation:

\begin{lem} \label{criterion}

Let $\C$ be a sheaf of categories on $\et(X)$, and $\A$ be a sheaf of abelian
groups on $\et(X)$. Suppose that for every
$(U\to X)\in\et(X)$ and every $C\in\C(U)$, we are given an isomorphism
$\Aut_{\C(U)}(C)\simeq \A(U)$ such that the following conditions hold:

\smallskip

\noindent{\em (0)}
There exists a covering $U\to X$ such that $\C(U)$ is non-empty.

\smallskip

\noindent{\em (1)}
If $C_1\to C_2$ is an isomorphism between two objects in $\C(U)$, then
the induced isomorphism $\Aut_{\C(U)}(C_1)\simeq \Aut_{\C(U)}(C_2)$ is
compatible with the identification of both sides with $\A(U)$.

\smallskip

\noindent{\em (2)}
If $f:U'\to U$ is a morphism in $\et(X)$ and $C\in\C(U)$, then the map
$$f^*_\C:\Aut_{\C(U)}(C)\to \Aut_{\C(U')}(f^*_\C(C))$$ is compatible with
the restriction map $\A(U)\to\A(U')$.

\smallskip

\noindent{\em (3)}
For any $U\in \et(X)$ and any two $C_1,C_2\in\C(U)$, there exist a
covering $f:U'\to U$ such that the objects $f^*_\C(C_1)$ and
$f^*_\C(C_2)$ of $\C(U')$ are isomorphic.

\smallskip

Then $\C$ has a canonical structure of a gerbe over $\Tors_\A$.
\end{lem}

\ssec{} We claim that $\Higgs_{\tilX}$ satisfies the conditions of this
lemma. Condition (0) is
a tautology: locally the cameral cover $\tilX\to X$
is induced from the universal one by means of a map $X\to\GNb$.

\medskip

Let $(E_G,\sigma,t)$ be an object of $\Higgs_{\tilX}(U)$. We must
construct an isomorphism
$$\on{Aut}_{\Higgs_{\tilX}(U)}(E_G,\sigma,t)\simeq T_{\tilX}(U).$$
Let us first assume that $E_G$ is trivialized and our Higgs bundle
corresponds to a map
$U\to\GNb$ such that $\tilU\overset{t}\simeq\GTb\underset{\GNb}\times U$.

In this case, an automorphism of $(E_G,\sigma)$ as an object of
$\Higgs(U)$ is the same as a map
$U\to \on{St}$ (cf. \secref{secstabilizers})
that covers the given map $X\to \GNb$. Now, \propref{stabilizer}
implies that this automorphism belongs to
$\on{Aut}_{\Higgs_{\tilX}(U)}(E_G,\sigma,t)$
if and only if the above map factors as $U\to \Z\to \on{St}$.

Now we apply \thmref{centralizers} which says
that $\Hom_{\GNb}(U,\Z)=\Hom_{\GNb}(U,\TT)=T_{\tilX}(U)$.

The fact that the map $\chi:\Z\to\TT$ is $G$-equivariant implies that our
isomorphism
between $\on{Aut}_{\Higgs_{\tilX}(U)}(E_G,\sigma,t)$ and $T_{\tilX}(U)$ is
independent
of the choice of a trivialization of $E_G$. In particular, by SC-1, it
defines the required
isomorphism for all $E_G$. The fact that conditions (1) and (2) are
satisfied is automatic from
the construction.

\medskip

Finally, let us check condition (3). Let $(E^1_G,\sigma^1,t^1)$ and
$(E^2_G,\sigma^2,t^2)$
be two objects of $\Higgs_{\tilU}(U)$. Without restricting the generality
we can assume
that both $E^1_G$ and $E^2_G$ are trivialized.

In this case, the data of $(\sigma^1,\sigma^2,t^1\circ (t^2)^{-1})$
defines a $U$-point of the scheme
$\H$. By \thmref{conjugation} we can locally find a map $g_U:U\to G$ which
conjugates $(\sigma^1,t^1)$ to
$(\sigma^2,t^2)$. We can regard $g_U$ as a gauge transformation, i.e. a
map $E_G^1\to E^2_G$, which defines
an isomorphism between $(E^1_G,\sigma^1,t^1)$ and $(E^2_G,\sigma^2,t^2)$.

\smallskip

Thus, \thmref{first} is proved.

\section{Proof of \thmref{descrofpoles}} \label{proofpoles}

\ssec{} It remains to prove \thmref{second}. In this section we will prove
\thmref{descrofpoles} which takes care of the universal situation.

\ssec{Step 1} First we show that our map $s_i^*(\L_{can})\to \L_{can}$
is an isomorphism off $D^{\alpha_i}$. To do that
let us analyze more closely the situation described in
\secref{Levicompat}.

\smallskip

Let $\Delta_J\subset\Delta$ be a root subsystem and let $M=M_J$ be the
corresponding
standard Levi subgroup. Let $\B_M$ denote the flag variety of $M$ and
$B_M=B\cap M$,
$U_M=U\cap M$.

\medskip

It is well-known that there exists a canonical closed embedding
$W_M\backslash W\times \B_M\to \B$:

A point $\bof'\in\B$ belongs to $w\times \B_M$, if and only if $\bof'$ is
in relative position
$w$ with respect to $P=P_J$ (this makes sense, as $P$-orbits in $\B$ are
parametrized exactly by
$W_M\backslash W$) and $\bof'\cap \mf$ is a Borel subalgebra in $M$.

\medskip

Consider the restriction of the canonical $T$-bundle $\L_{\B}$ to
$W_M\backslash W\times \B_M$. It is easy to
see that its further restriction to the connected component $1\times \B_M$
identifies with $\L_{\B_M}$.

Let $w\in W$ be a {\it minimal} representative of its coset in $W_M\backslash W$.
The action of  $w$
defines a map $1\times{\B}_M\to w^{-1}\times {\B}_M$.
Let us consider the pull-back
$w^*(\L_{\B}|_{w\times {\B}_M})$ as a $T$-bundle on $1\times {\B}_M={\B}_M$.
Let $\tilde{w}\in N$ be an element that projects to $w\in W$.

\smallskip

\begin{lem} \label{minimal}
We have a canonical $M$-equivariant isomorphism
$$w^*(\L_{{\B}}|_{w\times {\B}_M})\simeq \L_{{\B}_M}.$$
\end{lem}

\begin{proof}

Both $w^*(\L_{\B}|_{w\times {\B}_M})$ and $\L_{{\B}_M}$ are
$M$-equivariant $T$-bundles on ${\B}_M$. To prove that they
are isomorphic, we must show that the two homomorphisms
$B\cap M\to T$ corresponding to the base point $\bof\in{\B}_M$
coincide.

However, this follows from the fact that
$w^{-1}\times \bof=\on{Ad}_{\tilde{w}^{-1}}(\bof)$, which
is true since $w$ is minimal.

\end{proof}

\ssec{}

Let $w$ be as above. Consider the map
$$\overline{M/T}\overset{\text{\propref{pullbackGmodT}}}\longrightarrow
\GTb\overset{w}\to \GTb\to \B.$$
The fact that $\on{Ad}_{\tilde{w}^{-1}}(B)\cap M=B_M$ implies that the
above map coincides with
$$\overline{M/T}\to 1\times \B_M\overset{w}\to w\times \B_M\to \B.$$

Therefore, from \lemref{minimal} we obtain an isomorphism
$$\gamma'_{can}(\tilde{w}): w^*(\L_{can})|_{\overline{M/T}}\simeq
\L_{can}|_{\overline{M/T}}.$$

Moreover, it is easy to see that the above isomorphism
is induced by the restriction to $\overline{M/T}$
of the (meromorphic) isomorphism $\gamma_{can}(\tilde{w})$. In particular,
the {\it a priori} meromorphic isomorphism $\gamma_{can}(\tilde{w})$ is
regular
on $\overline{M/T}$.

\bigskip

Let us now go back to the situation of the theorem. We must check that the
meromorphic map
$s_i^*(\L_{can})\to \L_{can}$ has no poles along $D^{\alpha}$ if
$\alpha\neq\alpha_i$.
Choose a minimal Levi subgroup $M_j$ such that $w(\alpha_j)=\alpha$ for
some $w\in W$.

Then the fact that $\alpha\neq\alpha_i$ implies that
both $w$ and $s_i\cdot w$ are minimal representatives of
the corresponding cosets in $W/\langle s_j\rangle$. Then the above
discussion
shows that $s_i^*(\L_{can})\to \L_{can}$ has no poles on $w\times
\overline{M_j/T}$.

This proves what we need, since the $G$-orbit of $w\times
\overline{M_j/T}$ contains an open
part of $D^{\alpha}$ (cf. \propref{orbits}).

\ssec{Step 2}

Thus, we have shown that the poles of the map $s_i^*(\L_{can})\to
\L_{can}$
can occur only on $D^{\alpha_i}$. Let $M_i$ be the corresponding minimal
Levi
subgroup. As we have seen before, there is a natural embedding
$\overline{M_i/T}\to \GTb$, and $\L_{can}$ restricts to the corresponding
$T$-bundle on $\overline{M_i/T}$.

Since $s_i^*(\L_{can})\to \L_{can}$ is $G$-equivariant, to
determine the contribution of the divisor $D^{\alpha_i}$, it is enough
to perform the corresponding calculation for $M_i$. The latter
case easily reduces to $SL(2)$.

For $SL(2)$, $\GTb\simeq {\mathbb P}^1\times {\mathbb P}^1$
and $\L_{can}\simeq \O(1)\boxtimes \O$. Moreover, $-1\in S_2=W$
acts on $\GTb$ by swapping the two ${\mathbb P}^1$ factors, with
the fixed-point locus $\GTb^{-1}$ being the diagonal $\PP^1$.
Hence, $(-1)^*(\L_{can})\simeq \O\boxtimes \O(-1)$.

Therefore, we have a meromorphic map between $\O\boxtimes \O(-1)$ and
$\O(1)\boxtimes \O$, which is allowed to have zeroes and poles only
on the diagonal. Then it must have a zero of order $1$, by
degree considerations.

\section{Proof of \thmref{second}} \label{proof_second}

\ssec{The natural functor} Finally, we are ready to complete the proof
of the main result. First, we claim that there
is a natural functor $\Upsilon:\Higgs_{\tilX}\to \Higgs'_{\tilX}$:

Let $(E_G,\sigma,t)$ be an object of $\Higgs_{\tilX}(U)$,
where $\sigma:E_G\to \GNb$ is a $G$-equivariant map.
We can pull-back the universal object of $\Higgs'_{\GTb}(\GNb)$ (cf.
\secref{univex})
and obtain a $G$-equivariant object of
$\Higgs'_{\tilE_G}(E_G)$, where $\tilE_G$ is the induced cameral cover
of $E_G$.

By descent, it gives rise
to an object of $\Higgs'_{\tilX}(U)$ and this assignment is clearly a
functor between sheaves
of categories.

The key fact now is that $\Higgs'_{\tilX}$ is also a
gerbe over $\Tors_{T_{\tilX}}$. Condition (0)
of \lemref{criterion} follows from the mere existence of the functor
$\Upsilon$
and the fact that $\Higgs_{\tilX}$ satisfies condition (0).

Let $(\L,\gamma,\beta_i)$ be an object of $\Higgs'_{\tilX}(U)$. We must
identify the
group of its automorphisms with $T_{\tilX}(U)$.

By definition, this group consists of $T$-bundle automorphisms, which
respect the data of $\gamma$ and $\beta_i$.
However, a $T$-bundle map $\L\to\L$ is the same as a map $\tilU\to T$ and
compatibility with $\gamma$
implies that this map is $W$-equivariant. Therefore, we obtain a section
of $\overline{T}_{\tilX}(U)$.
Now, compatibility with $\beta_i$ is exactly condition (*).
(Recall that it suffices to impose condition (*)
for one representative in every $W$-orbit on the set of roots.
In particular, it is sufficient to impose it for simple roots only.)

It is easy to see that conditions (1) and (2) hold for the above
identification of
$\on{Aut}_{\Higgs'_{\tilX}(U)}(\L,\gamma,\beta_i)\simeq T_{\tilX}(U)$. In
addition, it follows
from the construction of $\chi$, that $\Upsilon:\Higgs_{\tilX}\to
\Higgs'_{\tilX}$ respects
the identifications of groups of automorphisms of objects with
$T_{\tilX}(U)$.

Assume for a moment that condition (3) of \lemref{criterion} has been
checked.
We claim, that this already implies \thmref{second}, because of the
following general fact:

\begin{lem}

Let $\C_1$ and $\C_2$ be two gerbes over $\Tors_\A$ and let
$\Upsilon:\C_1\to\C_2$ be a functor
between the corresponding sheaves of categories. Assume that for every
$U\in\et(X)$ and $C\in\C_1(U)$ we have a commutative square:
$$
\CD
\A(U) @>>> \on{Aut}_{\C_1(U)}(C) \\
@V{\on{id}}VV        @V{\Upsilon}VV   \\
\A(U) @>>> \on{Aut}_{\C_2(U)}(\Upsilon(C)).
\endCD
$$
Then $\Upsilon$ is an equivalence of $\Tors_\A$-gerbes.

\end{lem}

\ssec{The homogeneous version: $\Tors'_{T_{\tilX}}$} \label{homog}
It remains to prove that every two objects of
$\Higgs'_{\tilX}(U)$
are locally isomorphic. For that purpose we will introduce a sheaf of
Picard categories $\Tors'_{T_{\tilX}}$,
which will be the ``homogeneous'' version of $\Higgs'_{\tilX}(X)$.

Objects of $\Tors'_{T_{\tilX}}(U)$ are triples
$$(\L_0,\gamma_0,\beta_{i,0}),$$
where $(\L_0,\gamma_0)$ is a strongly $W$-equivariant $T$-bundle on
$\tilU$
and each $\beta_{i,0}$ is a trivialization of
$\alpha_i(\L_0)|_{D^{\alpha_i}_U}$.

The following compatibility conditions must hold:

\smallskip

\noindent (1) For a simple root $\alpha_i$, the data of
$\gamma_0(s_i):s^*_i(\L_0)\simeq \L_0$ defines,
after restriction to $D^{\alpha_i}_U$, a trivialization
$$\check\alpha_i(\alpha_i(\L_0)|_{D^{\alpha_i}_U})\simeq
\check\alpha_i(\O_{D^{\alpha_i}_U}).$$
We need that this trivialization coincides with
$\check\alpha_i(\beta_{i,0})$.

\smallskip

\noindent (2) Assume that $w\in W$ conjugates a simple root $\alpha_i$
to another simple root $\alpha_j$. The pull-back of $\beta_{j,0}$ under
$w$
is a trivialization of $\alpha_i(w^*(\L_0))|_{D^{\alpha_i}_U}$,
which via $\gamma_0(w)$
defines a trivialization of $\alpha_i(\L_0)$. Our condition is that
this trivialization coincides with $\beta_{i,0}$.

\smallskip

Morphisms in $\Tors'_{T_{\tilX}}(U)$ are by definition maps between
strongly $W$-equivariant $T$-bundles, compatible with the data of
$\beta_{i,0}$.

\medskip

If $(\L^1_0,\gamma^1_0,\beta^1_{i,0})$ and
$(\L^2_0,\gamma^2_0,\beta^2_{i,0})$ are two
objects of $\Tors'_{T_{\tilX}}(U)$ we can form their tensor product
$(\L^1_0\otimes \L^2_0,\gamma^1_0\otimes \gamma^2_0, \beta^1_{i,0}\otimes
\beta^2_{i,0})$
which will be a new object of $\Tors'_{T_{\tilX}}(U)$. Moreover, if
$(\L_0,\gamma_0,\beta_{i,0})$ is an object of $\Tors'_{T_{\tilX}}(U)$ and
$(\L,\gamma,\beta_i)$ is an object of $\Higgs'_{{\tilX}}(U)$, we can take
their tensor
product and obtain another object of $\Higgs'_{{\tilX}}(U)$.

It is easy to see that the above constructions define on
$\Tors'_{T_{\tilX}}$ a structure of a
sheaf of Picard categories and on $\Higgs'_{{\tilX}}$ a structure of a
gerbe over it.

Therefore, to prove that every two objects of $\Higgs'_{\tilX}(U)$
are locally isomorphic, it is enough to show that any object of
$\Tors'_{T_{\tilX}}(U)$
is locally isomorphic to the unit object, i.e. to the one with $\L_0$
being {\it the trivial}
$T$-bundle and $(\gamma_0,\beta_{i,0})$ being the tautological maps. The
last assertion
is equivalent to:

\begin{prop}  \label{TorsIsTors}
$\Tors'_{T_{\tilX}}$ is equivalent as a sheaf of Picard categories to
$\Tors_{T_{\tilX}}$.
\end{prop}
We proceed to prove this Proposition by showing that any object in
$\Tors'_{T_{\tilX}}(U)$ is locally isomorphic to the unit object.

\smallskip

\ssec{Step 1} Without restricting the generality, we can assume that $U=X$
and we must find an
\'etale covering $X'\to X$, over which a given object
$(\L_0,\gamma_0,\beta_{i,0})$
becomes isomorphic to the trivial one.

Fix a $\CC$-point $x\in X$. First, we will reduce our situation to the
case
when the ramification over $x$ is maximal
possible, i.e. when $x$ belongs to the image of $\underset{\alpha}\cap \,
D^\alpha_X$, where the intersection is taken over all roots of $G$.

After an \'etale localization we can assume that we have a map $X\to
\tf/W$
so that $\tilX=X\underset{\tf/W}\times \tf$.
Let $t$ be a point in $\tf$ which has the same image in $\tf/W$ as $x$.

By conjugating $t$, we can assume that there exists $J\subset I$ such that
$\alpha_j(t)=0$ for $j\in \Delta_J$ and
$\beta(t)\neq 0$ for $\beta\notin\Delta_J$.

We have a Cartesian square
$$
\CD
W\overset{W_J}\times(\tf\setminus \underset{\beta\notin
\Delta_J}\cup\,\tf^\beta) @>>>  \tf \\
@VVV        @VVV \\
(\tf\setminus \underset{\beta\notin
\Delta_J}\cup\,\tf^\beta)/W_J @>>> \tf/W,
\endCD
$$
in particular, the map $\tf/{W_J}\to \tf/W$ is \'etale in a neighbourhood
of the image of $t$
in $\tf/W$.

Therefore, the base change $X\Rightarrow X':=X\underset{\tf/W}\times
\tf/W_J$
is \'etale in a neighbourhood of $x$. This reduces us to the situation,
when
the $W$-cover $\tilX$ is induced from a $W_J$-cover $\tilX_J$, i.e.
$\tilX\simeq W\overset{W_J}\times\tilX_J$.

By restricting $(\L_0,\gamma_0,\beta_{i,0})$ to $\tilX_J$ we obtain an
object of
$\Tors'_{T_{\tilX_J}}(X)$ equivariant with respect to $W_J$. Moreover, it
is easy to see that this
establishes an equivalence between $\Tors'_{T_{\tilX}}$ and
$\Tors'_{T_{\tilX_J}}$,
thereby reducing us to the situation when $\Delta_J=\Delta$.

\ssec{Step 2}

According to Step 1, we may assume that there exists a unique
geometric point $\widetilde{x}\in\tilX$ over $x$.
To prove the assertion of the proposition, we can replace $X$ by
the spectrum of the local ring of $X$ at $x$. In this case all the
$D^{\alpha}_X$'s
and $\tilX$ are local too.

\smallskip

Let us choose a trivialization of our line bundle $\L_0$, subject only to
the
condition that it is compatible with the data of
$\beta_{i,0}$ at $\widetilde{x}$ for every simple root $\alpha_i$. We must
show
that this trivialization can be modified so that it
will be compatible with the structure on $\L_0$ of a $W$-equivariant
$T$-bundle, i.e. with the data of $\gamma_0$. (The argument given below
mimics
the proof of \propref{picardforgln}).

\smallskip

The discrepancy between our initial trivialization and $\gamma_0$
is given by a $1$-cocycle $\mu:W\to \Hom(\tilX,T)$.

The evaluation at $\widetilde{x}$ gives rise to a surjection of
$W$-modules:
$\Hom(\tilX,T)\to T$. Thus we obtain
a short exact sequence: $$0\to K\to \Hom(\tilX,T)\to T\to 0,$$
where $K$ consists of maps $\tilX\to T$ which have value $1$ at
$\widetilde{x}$.

Now, our condition on the trivialization (i.e. its compatibility with
$\beta_{i,0}$) and condition (1)
in the definition of $\Tors'_{T_{\tilX}}$ imply that $\mu(s_i)\in K$
for every simple reflection $s_i$. Hence, $\mu$ takes values in $K$.

However, since $\tilX$ is local, $K$ is torsion-free and divisible!
Hence, $H^1(W,K)=0$.

\medskip

Therefore, we can choose a trivialization of $\L_0$ which respects the
$W$-equivariant structure and the data of $\beta_{i,0}$ at
$\widetilde{x}$.
But this implies that it is compatible with the data of $\beta_{i,0}$ on
the entire $D^{\alpha_i}_X$, $\forall\,i\in I$.

Indeed, a possible discrepancy takes values in $\pm 1$, and its value is
constant along every connected component of $D^{\alpha_i}_X$. However,
by construction, each $D^{\alpha_i}_X$ is local with $\widetilde{x}$
being its unique closed point.

The proof of \propref{TorsIsTors}, and hence of \thmref{second}, is now
complete.

\ssec{Variant}
As was the case for $\on{Higgs}'$, we can give a much simplified description
of $\Tors'_{T_{\tilX}}$ in case our group $G$ does not have an
$\on{SO}(2n+1)$ direct factor.
In this case the data consists of a strongly equivariant $T$-bundle
$(L_0,\gamma_0)$, such that for a simple root $\alpha_i$
and a weight $\lambda$ orthogonal to the corresponding
coroot, the isomorphism
$\lambda(s_i^*(L))|_{D^{\alpha_i}}\to \lambda(L)|_{D^{\alpha_i}}$
induced by $\gamma(s_i)$ coincides with the tautological one.

\part{Some applications}
The point of our abstract notion of a Higgs bundle, as defined in
\secref{intrHiggs}, is that it provides a uniform approach to the analysis
of various more concrete objects. In the final sections we illustrate the
applications to Higgs bundles with values in a line bundle or in an
elliptic fibration.

\section{Higgs bundles with values}\label{values}

\ssec{}\label{c-and-c_X}
In \secref{intrHiggs} we defined a Higgs bundle over a scheme $X$ to be
a pair $(E_G,\sigma)$, where $\pi:E_G \to X$ is a
principal $G$-bundle over $X$, and $\sigma$ is a $G$-equivariant map
$\sigma:E_G\to \GNb$. We noted there that on a given $G$-bundle $E_G$,
a Higgs bundle is specified by a vector subbundle $\z_X$ of $\gf_{E_G}$
whose fibers are regular centralizers. (Recall that
$\gf_{E_G}:=E_G{\underset{G}\times}\gf$ is the adjoint bundle of $E_G$.)
In subsection \ref{centralizer_families} we defined the universal
centralizer $\z \subset \gf\times \GNb $, corresponding to the universal
Higgs bundle over $\GNb$. The family of centralizers $\z_X$ of
a general Higgs bundle $(E_G,\sigma)$ over $X$ is related to the
universal $\z$ by: $\pi^*\z_X = \sigma^*\z$, an equality of vector subbundles
of the trivial bundle $\gf \times E_G$ on the total space of $E_G$.
We recall also that by
Theorem \ref{Liealgebrasfirst}, $\z$ is isomorphic to $\tt = \on{Lie}(\TT)$.

Let $K$ be a line bundle on our base $X$. In the literature, the most
common notion of a Higgs bundle is:

\begin{defn}\label{KHiggs}
A $K$-valued Higgs bundle on $X$ is a pair $(E_G, s)$,
where $E_G$ is a principal $G$-bundle on $X$ and $s$ is a section
of $\gf_{E_G} \otimes K$.
\end{defn}

The section $s$ of $\gf_{E_G}\otimes K$ is called {\it regular} at a
point $x \in X$ if the corresponding local section of $\gf_{E_G}$
determined by some (hence, any) trivialization of $K$ at $x$ is regular.
We work instead with the following more general notion, which is
also better adapted to our setup.

\begin{defn}\label{RegKHiggs}
A regularized $K$-valued Higgs bundle on $X$ is a triple $(E_G,\sigma,s)$,
with $(E_G,\sigma)$ a Higgs bundle on $X$ and $s$ a section
of $\z_X\otimes K$, where $\z_X$ is the regular centralizer subbundle of
the adjoint bundle $\gf_{E_G}$ determined by $\sigma$.
\end{defn}

\ssec{Regular vs. regularized}\label{regular/ized}
A regularized $K$-valued Higgs bundle $(E_G,\sigma,s)$ on $X$ clearly
determines the unique $K$-valued Higgs bundle $(E_G,s)$ on $X$. Conversely,
if the section $s$ of $\gf_{E_G}\otimes K$ is everywhere regular, then we
can recover $\z_X\subset \gf_{E_G}$ as the centralizer of $s$, which
defines a regularized $K$--valued Higgs bundle. When $s$ is generically
regular, the family $\z_X$ of centralizers is still unique, if it exists.
In general, when $s$ is not necessarily regular, our definition adds to
the pair $(E_G,s)$ a choice of a regular centralizer containing $s$.

We want to establish the following result:

\begin{thm}   \label{Higgs with values}
A regularized $K$-valued Higgs bundle on $X$ is the same as a triple:

\smallskip

\noindent (a) A cameral cover $\tilX\to X$,

\smallskip

\noindent (b) A $W$-equivariant map $v:\tilX\to \tf\otimes K$ (of
schemes over $X$).

\smallskip

\noindent (c) An object of $\on{Higgs}'_{\tilX}(X)$,

\end{thm}

\begin{proof}

Given \thmref{second}, it remains to show that
the data (b) of a $W$-equivariant ``value''
map $\tilX\to \tf\otimes K$ is the same as the data
of a section $s$ of $\z_X \otimes K$. And indeed, giving such a section
$s: X \to \z_X \otimes K$
is equivalent to giving a $G$-equivariant section
${\tilde{s}}: E_G \to \sigma^*\z \otimes K$
of the pullback $\pi^*\z_X \otimes K= \sigma^*{\z} \otimes K$ over
$E_G$, cf. \ref{c-and-c_X} above.
By Theorem \ref{Liealgebrasfirst}, this is the same as a $G$-equivariant
section
${\tilde{s}'}: E_G \to \sigma^*\tt \otimes K$.
Now by the definition of
$\TT$ (cf. subsection \ref{TT}),
$\Hom_{\GNb}(E_G,\tt) = \Hom_W(\tilE_G, \tf)$.
Here $\tilE_G:=E_G \times_X\tilX$ is the $G$-equivariant cameral cover
of $E_G$ associated to the Higgs bundle on $E_G$ which is $\pi^*$ of
our given Higgs bundle $(E_G, \sigma)$ on $X$.
The section $\tilde{s}'$, and hence also our original section $s$, are
therefore equivalent to a $W$-equivariant map of $X$-schemes
$\overline{s}: \tilE_G \to \tf \otimes K$ which is also $G$-invariant.
But this is the same as a $W$-equivariant map of $X$-schemes
$v:\tilX \to \tf \otimes K$, as claimed.
\end{proof}

Note that in the data $(E_G,\sigma,s)$, the section
$s:X\to \z_X\otimes K$
is regular if and only if the corresponding map $v$
is an embedding. This follows from \lemref{regcriterion}. So we have:

\begin{cor}
A regular $K$-valued Higgs bundle on $X$ is the same as a triple:

\smallskip

\noindent (a) A cameral cover $\tilX\to X$,

\smallskip

\noindent (b) A $W$-equivariant embedding $v:\tilX\to \tf\otimes K$ (of
schemes over $X$).

\smallskip

\noindent (c) An object of $\on{Higgs}'_{\tilX}(X)$,

\end{cor}

\medskip

\ssec{The Hitchin map}
To conclude our discussion of $K$-valued Higgs bundles,
let us note that
the data (a) and (b) in the above theorem can be assembled into what can
be called ``a point of the Hitchin base''.

Assume that $X$ is proper, and let
${\mathbf B}(X,K)$ denote the algebraic stack which classifies
the data (a) and (b) of \thmref{Higgs with values}. I.e., for
a scheme $S$, $\Hom(S,{\mathbf B}(X,K))$ is the category
of pairs $(\tilX_S,v:\tilX_S\to \tof\otimes K)$, where
$\tilX_S$ is a cameral cover of $S\times X$, and $v$ is
a $W$-equivariant morphism of $X$-schemes.

On the other hand, let $\bfHiggs(X,K)$ denote the algebraic stack of
all regularized $K$-valued Higgs bundles on $X$. The {\it Hitchin map}
$h: \bfHiggs(X,K) \to {\mathbf B}(X,K)$ sends a regularized $K$-valued Higgs
bundle $(E_G,\sigma,s)$ given by data (a),(b) and (c) to the point of the
Hitchin base given by data (a) and (b).

\begin{cor}  \label{HitchFibers}
The fibers of the Hitchin map
$h:\bfHiggs(X,K)\to {\mathbf B}(X,K)$ can be identified (as categories) with
$\on{Higgs}'_{\tilX}(X)$.
By Corollary \ref{its_a_torsor}, the set of isomorphism classes of objects
of this fiber is a torsor over the abelian
group $H^1(X,T_{\tilX})$, and the torsor class is given in \thmref{second}.

\end{cor}

Note that our description of the fiber of the Hitchin map is independent
of the line bundle $K$.

\ssec{}

Let now $\bfHitch(X,K)$ denote the scheme
of sections of the fibration
$(\tf \otimes K)/W \to X$. In fact, $\bfHitch(X,K)$ is
non-canonically isomorphic to an affine space.

The relation between $\bfHitch(X,K)$ and ${\mathbf B}(X,K)$ is similar
in some respects to the relation between the vector space ${\tf} / W$
parametrizing semisimple adjoint orbits in the Lie algebra $\gf$ and the
stack $\gf / G$ of all $G$-orbits in $\gf$. In both cases, there is an open 
embedding of the variety into the stack, and there is a retraction of the 
stack onto the variety which is the identity on the variety.

In our case, the retraction $r:{\mathbf B}(X,K)\to \bfHitch(X,K)$ associates
to $v:\tilX\to \tf\otimes K$ the corresponding map $X\to (\tf \otimes K)/W$.
As for the open embedding $i:\bfHitch(X,K)\to {\mathbf B}(X,K)$:
starting with $X\to (\tf \otimes K)/W$, we recover $\tilX$ as
$$\tilX:=X\underset{(\tf \otimes K)/W}\times (\tf \otimes K),$$
and $v:\tilX\to \tf\otimes K$ is the second projection.

Obviously, the image $i(\bfHitch(X,K)) \subset {\mathbf B}(X,K)$ is the 
open substack corresponding to the condition that the map 
$\tilX\to \tof\otimes K$ is an embedding. By \corref{HitchFibers}, the 
preimage of $\bfHitch(X,K)\subset {\mathbf B}(X,K)$
under the Hitchin map is exactly the open substack of regular $K$-valued 
Higgs bundles. Let $R$ denote some regularized $K$-valued Higgs bundle on 
$X$. Note that the image $h(R) \in {\mathbf B}(X,K)$ determines whether 
$R$ is regular. A point in $\bfHitch(X,K)$, on the other hand, can be
the image of both regular and irregular $R$'s.

\ssec{Variant}
Definition \ref{RegKHiggs}, Theorem \ref{Higgs with values} and Corollary
\ref{HitchFibers}, remain unchanged
if we allow $K$ to be a vector bundle,
as is done in
\cite{Si} where $K=\Omega^1_X$ is the cotangent bundle. In Definition
\ref{KHiggs}, on the other hand, commutativity is not built in, so we
must impose it by hand: the components of the section $s$, with respect
to any local decomposition of $K$ as a sum of line bundles, must commute
with each other. Equivalently, the bracket of $s$ with itself,
interpreted as a section of $\gf_{E_G} \otimes \wedge^2 K$, must
vanish.

\section{Elliptic fibrations}\label{elliptics}

Let $f:Y \to X$ be a projective, flat, dominant morphism with integral
(that is, reduced and irreducible) fibers. Eventually we will specialize
this to the case of an elliptic fibration, but for now we will work with
the general situation.
We want to describe an application of our results to the study of regularized
$G$-bundles on $Y$ in terms of data on the base $X$ and along the (eventually,
elliptic) fibers.

By a regularization of a $G$-bundle $E_G$ on $Y$ we
mean a reduction of its structure group {\it along each fiber} to some regular
centralizer. In other words, we want a Higgs bundle $(E_G, \sigma)$ on
$Y$ whose group scheme of centralizers $\Z_Y$ (equivalently, its cameral cover
$\tilY \to Y$) is the pullback of some group scheme of centralizers
$\Z_X$ on $X$ (respectively, of a cameral cover $\tilX \to X$). More precisely:

\begin{defn}
A regularized $G$-bundle on $Y$ consists of the data $(\tilX,E_G,\sigma)$,
with $\tilX \to X$ a cameral cover of $X$, and
$(E_G,\sigma) \in \Higgs_{\tilY}(Y)$ a Higgs bundle
on $Y$  with cameral cover $\tilY :=\tilX\underset{X}\times Y$.
\end{defn}

In the case of an elliptic fibration, there is a natural notion of what it
means for a bundle (on $Y$) to be regular above a point (of $X$).
In analogy with the situation for $K$-valued Higgs bundles considered in
Subsection \ref{regular/ized}, ``most'' $G$-bundles on an elliptic curve are
indeed regular, and a regular bundle has a unique regularization. We review
these well-known facts below.

\ssec{} In general, our current situation is the analogue of Higgs
bundles with values,
in which we replace the bundle $K$ of values from
\secref{values} by the relative Picard scheme
$\Pic(Y/X)$. The tensor product $\tf \otimes_{\CC} K$ can be identified
with $\Lambda \otimes_{\ZZ} K$, so we take its analogue to be
$\Lambda \otimes_{\ZZ} \Pic(Y/X) =: \on{Bun}_T(Y/X).$
(Here $\Lambda$ is the lattice of coweights.)
Similarly, we will need the analogue of $\z_X \otimes K$. This is the sheaf
of groups $\Tors_{Y/X} := \Tors_{\Z_Y,Y/X}$, the sheafification of the
presheaf on $Y$ given by:

$\;\;\;\;\;\;\;\;\;$ $U \mapsto \{ \Z_Y$-torsors on $U$ modulo pullbacks of $\Z_X$-torsors$\}.$

\noindent (As above, $\Z_X$ is the group scheme of
regular centralizer subgroups with Lie algebra $\z_X$, and
$\Z_Y := f^*(\Z_X)$. ) In fancier language, we could think of
$\Tors_{Y/X}$ as a sheaf of Picard groupoids. But its objects have no
automorphisms, so we are dealing in fact with a sheaf of abelian groups.
In more detail:

\medskip

We introduce the sheaf of Picard categories $\Tors_{\Z_Y,Y/X}$ on $\et(X)$
as ``$\Z_Y$-torsors on $Y$ modulo pull-backs of $\Z_X$-torsors''. The
definition of $\Tors_{\Z_Y,Y/X}$ is as follows:

First, consider the presheaf of categories $\Tors^{\on{pre}}_{\Z_Y,Y/X}$, whose
objects over $U\to X$ are torsors over $U\underset{X}\times Y$
with respect to the sheaf $T_{\tilY}$. Morphisms between two such torsors
$\tau'$ and $\tau''$ are pairs $(\tau_X,\sigma)$, where $\tau_X$ is a $T_{\tilX}$-torsor
on $U$ and $\sigma$ is an isomorphism $\tau''\to \tau'\otimes f^*(\tau_X)$.
(Since $f:Y\to X$ is dominant, and thus
$\Gamma(U,T_{\tilX})\to \Gamma(U\underset{X}\times Y,T_{\tilY})$ is an
injection, it is easy to see that the morphisms defined this way form a
set and not just a category.)

The presheaf $\Tors^{\on{pre}}_{\Z_Y,Y/X}$ satisfies the first sheaf
axiom, but not the second one, i.e. not every descent data is automatically
effective. By applying the standard sheafification procedure, we obtain
from $\Tors^{\on{pre}}_{\Z_Y,Y/X}$ a sheaf of Picard categories, which we
denote by $\Tors_{\Z_Y,Y/X}$.

Note, however, that since the morphism $f:Y\to X$ is projective, objects
of $\Tors_{\Z_Y,Y/X}$ have no non-trivial automorphisms, because for
every $U$ as above, the map $\Gamma(U,T_{\tilX})\to \Gamma(U\underset{X}\times Y,T_{\tilY})$
is in fact an isomorphism. Hence,
$\Tors_{Y/X} := \Tors_{\Z_Y,Y/X}$ is in fact a sheaf of groups.

We need an explicit description of this sheaf:

\begin{lem}\label{descr quotient}
There is a canonical identification:

\noindent
$\Tors_{Y/X}(X) = \{v \in \on{Mor}_W(\tilX,\Bun_T(Y/X))|\
\alpha_i \circ v_{|D_X^{\alpha_i}}  =1 \in \Pic(Y/X),
\forall$ $\alpha_i \in I \}.$
\end{lem}

\noindent(As always, $I$ denotes the set of simple roots $\alpha_i.$)

\begin{proof}
We identify $\Tors_{Y/X}$ and $\Tors'_{Y/X}$ using Proposition
\ref{TorsIsTors}. There
is a natural map $\iota: \Tors'_{Y/X}\to \Mor_W(\tilX,\Bun_T(Y/X))$, sending
a $T$-bundle on $\tilY =\tilX \underset{X} \times Y$ to its classifying
morphism $v$. This map $\iota$ is clearly injective, and its image is
contained in the RHS.

We still have to prove the surjectivity of $\iota$, i.e. to show that a
morphism $v$ in the RHS satisfies the two
compatibility conditions between $\beta$'s and $\gamma$'s stated in
\ref{homog}. It suffices to do so locally,
and then we may assume that $f:Y \to X$ has a section.
In this case, we can identify $\Tors'_{Y/X}$ with the sheaf of $T$-bundles
on $\tilY$ satisfying the two compatibility conditions between $\beta$'s
and $\gamma$'s and which additionally are
trivialized along the section $X \subset Y$. Similarly, we can identify
$\Bun_T(Y/X))$ with $T$-bundles on $Y$ which are trivialized along the
section.

Each of the compatibility conditions requires the
equality of two given trivializations of some ($T$- or ${\GG}_m$-) bundle
over $D_X^i\underset{X}\times Y$.
Now our assumption, $\alpha_i \circ v_{|D_X^i} = 1$, together with the
assumed trivialization of all objects along the section, guarantees that
these equalities hold over the section. The difference between the
two trivializations is therefore a global automorphism which equals the
identity along the section, so it is the identity everywhere since the
fibers of $f$ are integral and proper.

\end{proof}

\ssec{} \label{ses}
By construction, we have a short exact sequence of Picard categories:

$$0 \to \Tors_{T_{\tilX}} \to f_*(\Tors_{T_{\tilY}}) \to \Tors_{Y/X} \to 0.$$

As in Subsection \ref{Q}, an element $v \in \Tors_{Y/X}(X)$ determines a
$\Tors_{T_{\tilX}}$-gerbe which we denote $\C_v$.
In fact, for $(U\to X)\in \et(X)$, $\C_v(U)$ is the category of all possible
lifts of $v$ to a $T_{\tilY}$--torsor on $U\underset{X}\times Y$.

The main result of this section
is the following analogue of \thmref{Higgs with values}:

\begin{thm}   \label{bundles on elliptics}
A regularized $G$-bundle on $Y$ is the same as a triple:

\smallskip

\noindent (a) A cameral cover $\tilX\to X$,

\smallskip

\noindent (b) A $W$-equivariant map $v:\tilX\to \Bun_T(Y/X)$ (of $X$-schemes),
satisfying:

\noindent
$\alpha_i \circ v_{|D^{\alpha_i}}  =1 \in \Pic(Y/X), \forall$ simple root
$\alpha_i$, and

\smallskip

\noindent (c) An object of
$\on{Higgs}'_{\tilX}(X) \underset{\Tors_{T_{\tilX}}}\otimes \C_v$.

\end{thm}

\begin{proof}
Let us fix a cameral cover $\tilX\to X$ and
consider regularized $G$-bundles on $Y$
corresponding to this fixed $\tilX$ as
a sheaf of categories over $X$, denoted by $\on{Reg}_{\tilX}(Y)$.

By \thmref{first}, $\on{Reg}_{\tilX}(Y)$
is a gerbe over the sheaf of Picard categories
$f_*(\Tors_{T_{\tilY}})$. This gerbe is induced from
the $\Tors_{T_{\tilX}}$-gerbe $\Higgs_{\tilX}$
by the homomorphism
$\Tors_{T_{\tilX}} \to f_*(\Tors_{T_{\tilY}})$, cf. \secref{exact sequences}.

Thus, according to \lemref{induced gerbes}, we have a functor
$\on{Reg}_{\tilX}(Y)\to \Tors_{Y/X}$, and for a given
object $v\in \Tors_{Y/X}(X)$ the category-fiber of the
above functor is a $\Tors_{T_{\tilX}}$-gerbe, which can be canonically
identified with
$\on{Higgs}_{\tilX}(X) \underset{\Tors_{T_{\tilX}}}\otimes \C_v$.

Finally, according to \lemref{descr quotient},
an object $v\in \Tors_{Y/X}(X)$ is equivalent to data (b) above.

\end{proof}

\ssec{}

Now let us assume that $X$ is projective as well.
As our analogue of $\bfHiggs(X,K)$,
we will consider the algebraic stack $\bfReg(X,Y)$ which associates
to a scheme $S$ the category of regularized $G$-bundles on $S\underset{X}\times Y$
(with respect to the projection $S\underset{X}\times Y\to S$).

We can now describe an analogue of the Hitchin map. Indeed, let
${\mathbf B}(X,Y)$ be the stack whose $S$-points are pairs:
$(\tilX_S,v)$ consisting of a cameral cover of 
$S\times X$ and a $W$-equivariant map $v: \tilX_S\to \Bun_T(Y/X)$ 
of $X$-schemes.

We have a natural map of stacks
$h:\bfReg(X,Y)\to {\mathbf B}(X,Y)$.

\begin{cor}
The fiber of the spectral map $\bfReg(X,Y) \to {\mathbf B}(X,Y)$ over a cameral
point $(\tilX,v) \in {\mathbf B}(X,Y)$ can be identified with
$\on{Higgs}'_{\tilX}(X) \underset{\Tors_{T_{\tilX}}}\otimes \C_v$.
The set of isomorphism classes of objects
of this fiber is a torsor over the abelian
group $H^1(X,T_{\tilX})$.
\end{cor}

In the case of $K$-valued Higgs bundles, we saw in Corollary
\ref{HitchFibers} that the fiber of the Hitchin map
$\bfHiggs(X,K) \to {\mathbf B}(X,K)$  is independent
of the line bundle $K$.
Note in contrast that the fiber
$\on{Higgs}'_{\tilX}(X) \underset{\Tors_{T_{\tilX}}}\otimes \C_v$
of the spectral map could depend on the original map $f:Y \to X$. This
dependence is mild though: it affects only
the second factor, $\C_v$. A simplification occurs when $f:Y \to
X$ has a global section: in this case $\C_v$ is always trivial,
because its defining short exact sequence of Picard categories
\ref{ses} is split. It follows that the category
$\on{Reg}_{\tilX}(Y)$
of regularized bundles with a specified cameral cover $\tilX$ factors:
$$\on{Reg}_{\tilX}(Y) = \Tors_{Y/X} \times \on{Higgs}'_{\tilX}(X).$$

\ssec{}

In addition to the stack ${\mathbf B}(X,Y)$, one can also define
an analogue of the space $\bfHitch(X,K)$:
we let the space $\bfHitch(X,Y)$ denote the scheme of all sections of 
the fibration $(\Bun_T(Y/X))/W \to X$.
As before, we have an obvious retraction ${\mathbf B}(X,Y)\to \bfHitch(X,Y)$.
The analogue of the embedding $\bfHitch(X,K)\to {\mathbf B}(X,K)$
can be described as follows:

Consider the $W$-cover $\Bun_T(Y/X)\to(\Bun_T(Y/X))/W$
and let $({\Bun}_T(Y/X))^0/W$ be the maximal
open subscheme over which this cover is cameral; let
$\Bun_T(Y/X)^0$ denote its preimage in $\Bun_T(Y/X)$.

We will have to shrink $({\Bun}_T(Y/X))^0/W$ to a still smaller
open subscheme:

For a simple root $\alpha_i$ consider the corresponding ramification
divisor $$D^{\alpha_i}_{({\Bun}_T(Y/X))^0/W}\subset \Bun_T(Y/X)^0.$$
Under the map $\Bun_T(Y/X)\to \Pic(Y/X)$ given by $\alpha_i$,
the image of $D^{\alpha_i}_{({\Bun}_T(Y/X))^0/W}$ is contained
in the set of $\ZZ_2$-torsion points of $\Pic(Y/X)$.

We define the open subscheme $({\Bun}_T(Y/X))^{00}/W$ of
$({\Bun}_T(Y/X))^0/W$ by removing those points, whose
preimage in $\Bun_T(Y/X)^0$ maps to a {\it non-unit}
point in $\Pic(Y/X)$ by means of the above map.
Let ${\Bun}_T(Y/X))^{00}\to ({\Bun}_T(Y/X))^{00}/W$
denote the corresponding cameral cover.

\medskip

Finally, let
$\bfHitch(X,Y)^0$ be the open subscheme of
$\bfHitch(X,Y)$, which corresponds to sections
whose values belong to $({\Bun}_T(Y/X))^{00}/W$.

The fiber product construction gives the desired map
$i:\bfHitch(X,Y)^0\subset {\mathbf B}(X,Y)$. Its image
is the open substack corresponding to the locus
where the map $v:\tilX\to {\Bun}_T(Y/X)$ is
an embedding.

\ssec{The case of an elliptic fiber}

The main relevance of the above results is to the case that $f:Y \to X$
is an elliptic fibration. This is due to the existence in this case of a
good notion of a regular bundle, analogous to the notion of a regular
$K$-valued Higgs bundle. Take the group $G$ to be semisimple, and consider
the case of a single elliptic curve $Z$.

For any semistable $G$-bundle $E_G$ on $Z$, the dimension of the group
$H:=\on{Aut}_G(E_G)$ of (global) automorphisms of $E_{G}$ is $\geq r$. We
say that $E_G$ is regular if $\dim(H)=r$. In this case, $H$
is commutative and there exists an embedding $H\to G$
and a principal $H$-bundle $E_H$ on $Z$ such that
$E_G\simeq G\overset{H}\times E_H$. A regular bundle has a unique
regularization.

These results can be found in \cite{D3,FM,FMW2,L} and elsewhere.
In fact, the moduli space $M_G(Z)$ of  (S-equivalence classes of)
semistable, topologically trivial
$G$-bundles on the elliptic curve $Z$ is well understood. As a complex
variety, it is isomorphic to $M_T(Z)/W$.  (This is proved analytically
(e.g. \cite{FMW2}) using Borel's result that in a simply connected compact
group, any two commuting elements are contained in a maximal torus. An
algebraic proof was given in \cite{L}.) Each
S-equivalence class contains a unique regular representative as well as a
unique semisimple representative (i.e. one whose structure group can be
reduced to $T$). For a generic point of the moduli space, the
S-equivalence class consists of a unique isomorphism class, which is both
regular and semisimple. A similar but somewhat more complicated description
exists for all reductive $G$, cf. \cite{FM}. Returning to an elliptic
family $f:Y\to X$, we find ourselves in a situation analogous to that
which we had for $K$-valued Higgs bundles: a ``generic'' $G$-bundle on
$Y$ which is semistable along the elliptic fibers should be regular on
the generic fiber, and therefore its restriction to a dense open
$X_0 \subset X$ should admit a unique regularization to which we can
apply our results.

\end{document}